\documentclass[11pt]{amsart}
\usepackage[T1]{fontenc}
\usepackage{lmodern,xcolor}
\usepackage[letterpaper,textwidth=5.35in,textheight=8.05in,centering]{geometry}
\usepackage{amsmath,amssymb,amsthm,mathtools}
\usepackage{tikz}
\usetikzlibrary{arrows.meta}
\usepackage[hidelinks]{hyperref}
\hypersetup{pdftitle={The Adamczak--Prochno--Strzelecka--Strzelecki conjectures},pdfauthor={Daniel M. Pellegrino and Eduardo V. Teixeira},pdfsubject={Structured random matrices},pdfkeywords={Gaussian matrices, structured random matrices, operator norms, sequence spaces, variance profiles}}
\newtheorem{theorem}{Theorem}[section]
\newtheorem{lemma}[theorem]{Lemma}
\newtheorem{corollary}[theorem]{Corollary}
\newtheorem{proposition}[theorem]{Proposition}
\newcommand{\supp}{\operatorname{supp}}
\newcommand{\R}{\mathbb R}
\newcommand{\prob}{\mathbb P}
\newcommand{\norm}[1]{\left\|#1\right\|}

\numberwithin{equation}{section}
\title[The Adamczak--Prochno--Strzelecka--Strzelecki conjectures]{The Adamczak--Prochno--Strzelecka--Strzelecki conjectures}
\date{}
\author[D. M. Pellegrino]{Daniel M. Pellegrino}
\address{Departamento de Matem\'atica, Universidade Federal da Para\'iba, 58051-900 Jo\~ao Pessoa, PB, Brazil}
\email{daniel.pellegrino@academico.ufpb.br}
\author[E. V. Teixeira]{Eduardo V. Teixeira}
\address{Department of Mathematics, Oklahoma State University, Stillwater, OK, USA}
\email{eduardo.teixeira@okstate.edu}
\begin{document}
\begin{abstract}
We settle the two conjectures of Adamczak, Prochno, Strzelecka and Strzelecki concerning the expected norms of finite Gaussian matrices and the almost sure boundedness of infinite Gaussian matrices acting from $\ell_p$ to $\ell_q$, for all $1\le p,q\le\infty$.
\end{abstract}
\maketitle
\section{Introduction}
The use of random signs in norm estimates for trigonometric series goes back to Salem and Zygmund \cite{SZ}; see also Kahane \cite{Kahane}. For operators between sequence spaces, Bennett, Goodman and Newman \cite{BGN} studied the corresponding problem for random matrices. In 1975 they proved that an $m\times n$ matrix with independent centered entries bounded in modulus by one has expected $\ell_2^n\to\ell_q^m$ norm at most $C_q\max\{\sqrt n,m^{1/q}\}$ for $2\le q<\infty$. Independent random signs attain this order. In 1977, Bennett \cite{Bennett} obtained the optimal asymptotic order of the smallest possible norms of sign matrices acting from $\ell_p^n$ to $\ell_q^m$ for arbitrary $1\le p,q\le\infty$. Polynomial and multilinear versions of the Kahane--Salem--Zygmund inequality were studied in \cite{AR,MS}. Asymptotic estimates for the constants in the Kahane--Salem--Zygmund and Bennett inequalities were obtained in \cite{PR}.

For Gaussian matrices with independent entries of arbitrary variances, the problem is to estimate the mean norm in terms of the matrix of variances.

For the spectral norm, Seginer \cite{Seginer} proved in 2000 that, for independent and identically distributed centered entries, the mean norm is comparable to the sum of the expected maxima of the Euclidean row and column norms. Whether such a description remains valid for independent centered Gaussian entries with arbitrary variances was a conjecture of Lata\l a. His estimate in \cite{Latala2005} involved the largest square roots of the row and column sums of variances, as well as a term depending on the fourth moments of all the entries. Further progress on the spectral problem was made by Bandeira and van Handel \cite{BvH} and by van Handel \cite{vH}. The conjecture was eventually proved in 2018 by Lata\l a, van Handel and Youssef \cite{LVHY}, who also characterized the infinite Gaussian matrices defining bounded operators on $\ell_2$.

Meanwhile, Gu\'edon, Hinrichs, Litvak and Prochno \cite[Conjecture~1.2]{GHLP} had posed the analogous question for operators from $\ell_p$ to $\ell_q$ with $p\le2\le q<\infty$. Their 2017 paper contained upper estimates with a logarithmic loss and proposed a dimension-free estimate in terms of the largest appropriate row and column norms and the expected modulus of the largest entry. Outside this range of exponents, however, the same expression need not give the correct order. In \cite{APSS}, first circulated in 2021 and published in the 2024 volume of \emph{Mathematische Annalen}, Adamczak, Prochno, Strzelecka and Strzelecki identified the deterministic quantities appropriate to the other ranges of exponents and formulated two conjectures: one for the expected norm of finite Gaussian matrices and the other for the almost sure boundedness of infinite Gaussian matrices.

Let $m,n\in\mathbb N$, $1\le p,q\le\infty$ and $A=(a_{ij})_{i\le m,j\le n}\in\R^{m\times n}$. We use $1/\infty=0$, $\mathbb N=\{1,2,\ldots\}$ and the natural logarithm $\log$. For a vector $x$, $\supp x=\{j:x_j\ne0\}$. For a finite index set $I$, let $\R^I$ be the space of real families indexed by $I$. For $N\in\mathbb N$ we identify $\R^N$ with $\R^{\{1,\ldots,N\}}$. For $0<t<\infty$, write $\ell_t^N$ for $\R^N$ equipped with $\norm{x}_t=(\sum_{j=1}^N|x_j|^t)^{1/t}$; for $t=\infty$ use the supremum norm. When $t<1$, this is a quasi-norm. For $t\in(0,\infty]$ write $\ell_t^I$ for $\R^I$ with the corresponding (quasi-)norm and $B_t^I=\{x\in\R^I:\norm{x}_t\le1\}$; if $I=\{1,\ldots,N\}$, write $B_t^N$. For $J\subseteq I$, $x_J=(x_j)_{j\in J}\in\R^J$ denotes the restriction of $x\in\R^I$; when regarded as a vector in $\R^I$, it is extended by zero. For $a,b\in(0,\infty]$ and $V=(v_{ij})\in\R^{m\times n}$, the linear map and its norm are
\[
 V:\ell_a^n\longrightarrow\ell_b^m,\qquad
 Vx=\left(\sum_{j=1}^n v_{ij}x_j\right)_{i=1}^m,
 \qquad \norm{V}_{a\to b}=\sup_{x\in B_a^n}\norm{Vx}_b.
\]
For $I\subseteq\{1,\ldots,m\}$ and $J\subseteq\{1,\ldots,n\}$, the restriction $(v_{ij})_{i\in I,j\in J}$ defines a map $\ell_a^J\to\ell_b^I$. The space indexed by the empty set is the zero space. For $t\in[1,\infty]$, the conjugate exponent $t'$ is determined by $1/t+1/t'=1$.

For fixed $1\le p,q\le\infty$, define $D_1,D_2:\R^{m\times n}\to[0,\infty)$ by
\begin{equation}\label{d1}
 D_1(A)=\sup_{x\in B_p^n}\left[\sum_{i=1}^m
 \left(\sum_{j=1}^n a_{ij}^2 x_j^2\right)^{q/2}\right]^{1/q},
\end{equation}
and
\begin{equation}\label{d2}
 D_2(A)=\sup_{y\in B_{q'}^m}\left[\sum_{j=1}^n
 \left(\sum_{i=1}^m a_{ij}^2 y_i^2\right)^{p'/2}\right]^{1/p'}.
\end{equation}
Expressions such as \eqref{d1} and \eqref{d2} are interpreted as maxima if the outer exponent is infinite. The dependence of $D_1,D_2$ on $p,q$ is suppressed in the notation. For matrices $A,B$ of the same size write $A\circ B=(a_{ij}b_{ij})$ for the entrywise product; write $A^t$ for the transpose. Then
\begin{align}
 D_1(A)&=\norm{A\circ A}_{p/2\to q/2}^{1/2},\label{sq1}\\
 D_2(A)&=\norm{(A\circ A)^t}_{q'/2\to p'/2}^{1/2}.\label{sq2}
\end{align}
Indeed, $\norm{(x_j^2)_j}_{p/2}=\norm{x}_p^2$, and every nonnegative element of $B_{p/2}^n$ is the coordinatewise square of an element of $B_p^n$. Since $A\circ A$ has nonnegative entries, its operator norm can be computed on nonnegative vectors, by replacing any vector by its coordinatewise modulus. The same identity with $A^t$ and the exponents $(q',p')$ gives \eqref{sq2}. 

On a probability space $(\Omega,\mathcal A,\prob)$, let $g_{ij}:\Omega\to\R$, $1\le i\le m$, $1\le j\le n$, be independent random variables with standard Gaussian distribution (mean zero, variance one). We consider the random operator
\[
 G_A:\ell_p^n\longrightarrow\ell_q^m,\qquad
 G_Ax=\left(\sum_{j=1}^n a_{ij}g_{ij}x_j\right)_{i=1}^m.
\]
In particular, $a_{ij}^2$ is the variance of the $(i,j)$-th entry of $G_A$. 
For $t\in(0,\infty)$, let $L_t(\Omega,\mathcal A,\prob)$ be the space of real measurable functions $Z:\Omega\to\R$ with $\norm{Z}_{L_t}=(\int_\Omega|Z|^t\,d\prob)^{1/t}<\infty$, identified modulo equality almost everywhere. In particular, $\gamma_t=\norm{g}_{L_t}$ for a standard Gaussian $g$. The problem is to estimate
\begin{equation}\label{meannorm}
 \int_\Omega\norm{G_A}_{p\to q}\,d\prob
 =\int_\Omega\sup_{x\in B_p^n}\norm{G_Ax}_q\,d\prob
\end{equation}
in terms of the deterministic matrix $A$, uniformly in $m$ and $n$.

A first lower bound is obtained by fixing $x\in B_p^n$ before taking the integral. By independence, the random variable $(G_Ax)_i$ has the same distribution as $\sigma_i(x)g$, where
\[
 \sigma_i(x)^2=\int_\Omega|(G_Ax)_i|^2\,d\prob
 =\sum_{j=1}^n a_{ij}^2x_j^2.
\]
It follows that, for $q<\infty$,
\begin{equation}\label{fixedx}
 \left(\int_\Omega\norm{G_Ax}_q^q\,d\prob\right)^{1/q}
 =\gamma_q\left[\sum_i\left(\sum_j a_{ij}^2x_j^2\right)^{q/2}\right]^{1/q}.
\end{equation}
Moreover, $\int_\Omega|(G_Ax)_i|\,d\prob=\gamma_1\sigma_i(x)$. By convexity of the $\ell_q$-norm, also for $q=\infty$,
\begin{align*}
 \int_\Omega\norm{G_A}_{p\to q}\,d\prob
 &\ge\int_\Omega\norm{G_Ax}_q\,d\prob\\
 &\ge\norm{\left(\int_\Omega|(G_Ax)_i|\,d\prob\right)_{i=1}^m}_q\\
 &=\gamma_1\norm{(\sigma_i(x))_{i=1}^m}_q.
\end{align*}
Taking the supremum over $x\in B_p^n$ gives the lower bound $\gamma_1D_1(A)$. By duality, $\norm{G_A}_{p\to q}=\norm{G_A^t}_{q'\to p'}$. Applying the same calculation to the transpose gives the lower bound $\gamma_1D_2(A)$, and therefore
\begin{equation}\label{intro:basiclower}
 \int_\Omega\norm{G_A}_{p\to q}\,d\prob
 \ge\gamma_1\max\{D_1(A),D_2(A)\}
 \ge\frac{\gamma_1}{2}(D_1(A)+D_2(A)).
\end{equation}
The exponents in \eqref{sq1} and \eqref{sq2} thus come from the variance of a Gaussian sum and ordinary duality.

If $1\le p\le2\le q\le\infty$, the quantities in \eqref{sq1} and \eqref{sq2} simplify to
\begin{equation}\label{intro:crossing}
 D_1(A)=\max_j\norm{(a_{ij})_{i=1}^m}_q,\qquad
 D_2(A)=\max_i\norm{(a_{ij})_{j=1}^n}_{p'}.
\end{equation}
Indeed, for $x\in B_p^n$ we have $\sum_j x_j^2\le1$. Since $q/2\ge1$,
\[
 \norm{(\sum_j a_{ij}^2x_j^2)_i}_{q/2}
 \le\sum_j x_j^2\norm{(a_{ij}^2)_i}_{q/2}
 \le\max_j\norm{(a_{ij})_i}_q^2.
\]
Equality is attained at a suitable coordinate vector; the formula for $D_2$ follows by transposition. Thus the estimate conjectured in \cite{GHLP} for $q<\infty$ is the equivalence, up to constants depending only on $p$ and $q$, between \eqref{meannorm} and
\begin{equation}\label{intro:ghlp}
 D_1(A)+D_2(A)+\int_\Omega\max_{i,j}|a_{ij}g_{ij}|\,d\prob.
\end{equation}
In the spectral case $p=q=2$, this is precisely the estimate established in \cite{LVHY}. For general exponents the two square-function quantities need not reduce to maxima of row and column norms. For example, if $p=2$ and $1<q<2$, then
\begin{equation}\label{p2example}
 D_2(A)=\left[\sum_i\left(\sum_j a_{ij}^2\right)^{s/2}\right]^{1/s},
 \qquad s=\frac{2q}{2-q}.
\end{equation}
In fact, with $c_i=(\sum_j a_{ij}^2)^{1/2}$, we have
\[
 D_2(A)^2=\sup_{\norm{y}_{q'}\le1}\sum_i c_i^2y_i^2.
\]
Since $1/2=1/s+1/q'$, H\"older's inequality gives $D_2(A)\le\norm{c}_s$. If $c\ne0$, equality is attained at $y_i=(c_i/\norm{c}_s)^{s/q'}$. In particular, $D_2(A)$ can depend on the entire sequence of row norms, not just the largest one.

The lower bound \eqref{intro:basiclower} is not always of the correct order. Consider $A=I_n$, the identity matrix of order $n$. The operator $G_{I_n}$ is diagonal with independent standard Gaussian entries $g_{jj}$. If $p\le q$, then
\[
 \norm{G_{I_n}}_{p\to q}=\max_{1\le j\le n}|g_{jj}|,
 \qquad D_1(I_n)=D_2(I_n)=1.
\]
The equality for the random operator follows from $\norm{x}_q\le\norm{x}_p$ and evaluation at a unit vector. Since the mean of the maximum of $n$ independent standard Gaussian moduli is comparable to $\sqrt{\log(n+1)}$, no uniform upper bound by $D_1+D_2$ is possible in this range. In contrast, if $q<p$, put $s=(1/q-1/p)^{-1}$. H\"older's inequality, with equality at a suitable vector, gives
\[
 \norm{G_{I_n}}_{p\to q}=\left(\sum_{j=1}^n|g_{jj}|^s\right)^{1/s},
 \qquad D_1(I_n)=D_2(I_n)=n^{1/s}.
\]
The mean of the first quantity lies between $\gamma_1 n^{1/s}$ and $\gamma_s n^{1/s}$, by convexity of the $\ell_s$-norm and Jensen's inequality. 

For $q<2$ and $j\in\{1,\ldots,n\}$, the column $G_Ae_j$ satisfies
\begin{equation}\label{intro:colscale}
 \sup_{\norm{y}_{q'}\le1}\left(\int_\Omega
 |\langle G_Ae_j,y\rangle|^2\,d\prob\right)^{1/2}
 =\sup_{\norm{y}_{q'}\le1}\left(\sum_i a_{ij}^2y_i^2\right)^{1/2}
 =\norm{(a_{ij})_i}_{2q/(2-q)}=:b_j.
\end{equation}
The last identity follows from H\"older's inequality, since $1/2=1/q'+(2-q)/(2q)$. For each $j$ choose a vector $y^{(j)}\in B_{q'}^m$ attaining the supremum in \eqref{intro:colscale}. The random variables $\langle G_Ae_j,y^{(j)}\rangle$ are independent centered Gaussians of standard deviations $b_j$. Consequently,
\[
 \int_\Omega\norm{G_A}_{p\to q}\,d\prob
 \ge\int_\Omega\max_j\norm{G_Ae_j}_q\,d\prob
 \ge\int_\Omega\max_j b_j|g_j|\,d\prob,
\]
where $g_j:\Omega\to\R$, $1\le j\le n$, are independent standard Gaussians (on a product extension if necessary). Let $b^*$ denote the nonincreasing rearrangement of $(b_j)_j$. For each $1\le k\le n$, the last integral is at least $b_k^*\int\max_{1\le j\le k}|g_j|\,d\prob$. The Gaussian maximum estimate (see also \cite[Proposition 5.3]{APSS}) gives the lower bound, up to an absolute constant,
\[
 \max_{1\le j\le n}b_j^*\sqrt{\log(j+1)}.
\]
The dual calculation for $p>2$ gives a corresponding quantity involving the rows. Finally, since $\norm{G_A}_{p\to q}\ge\max_{i,j}|a_{ij}g_{ij}|$, the mean of the largest entry is always a lower bound.

Adamczak, Prochno, Strzelecka and Strzelecki showed that the expected largest entry in \eqref{intro:ghlp} does not always account for the fluctuations when $p\le q<2$ or $2<p\le q$ (see \cite[Section 5.3]{APSS}). They also showed that the column and row fluctuation terms are controlled by $D_1+D_2$ when $q<p$ (see \cite[Section 5.2]{APSS}). These observations, along with their sharp estimates at the boundary exponents, led to the following conjecture. Throughout, $U\asymp_{p,q} W$ means that $C_{p,q}^{-1}W\le U\le C_{p,q}W$ for some $C_{p,q}\ge1$ independent of $m,n,A$. Constants with subscripts may vary from line to line and depend only on the indicated parameters.

\medskip
\noindent\textit{Conjecture 1 of \cite{APSS}.} For $1\le p,q\le\infty$,
\begin{equation}\label{conj}
 \int_\Omega\norm{G_A}_{p\to q}\,d\prob
 \asymp_{p,q} D_1(A)+D_2(A)+L_{p,q}(A),
\end{equation}
where
\begin{equation}\label{lastterm}
 L_{p,q}(A)=\begin{cases}
 \displaystyle\int_\Omega\max_{i,j}|a_{ij}g_{ij}|\,d\prob,&p\le2\le q,\\[4pt]
 \displaystyle\max_{1\le j\le n} b_j^*\sqrt{\log(j+1)},&p\le q<2,\\[4pt]
 \displaystyle\max_{1\le i\le m} d_i^*\sqrt{\log(i+1)},&2<p\le q,\\[4pt]
 0,&q<p.
 \end{cases}
\end{equation}
Here $b_j$ is as in \eqref{intro:colscale}, $d_i=\norm{(a_{ij})_{j=1}^n}_{2p/(p-2)}$ for $p>2$ (with exponent $2$ at $p=\infty$), and $u^*$ denotes the nonincreasing rearrangement of any finite nonnegative vector $u$.
When $p\le q=2$, set $b_j=\max_i|a_{ij}|$. By \cite[Proposition~5.5]{APSS}, the largest-entry and ordered-column expressions are equivalent after $D_1(A)$ is added to each. Transposition gives the corresponding equivalence with $D_2(A)$ when $2=p\le q$.
\medskip

The lower estimate in \eqref{conj} was proved in \cite{APSS}. The authors obtained the reverse inequality when either $p$ or $q$ is $1$ or $\infty$; their general upper bounds for the interior exponents involved logarithmic factors in the dimensions. In 2026, Lata\l a and Strzelecka \cite[Theorem 2]{LS} established the conjecture of Gu\'edon, Hinrichs, Litvak and Prochno throughout the range $p\le2\le q$ (the boundary cases $p=1$ and $q=\infty$ were already known). A further question in \cite{LS} concerns the sharp dependence on $p,q$ in this range; a candidate proof of that quantitative estimate was subsequently proposed by Bednorz, Martynek and Meller \cite{BMM}. The estimates in the following two theorems concern $q<p$ and the two dual regions $1<p\le q<2$ and $2<p\le q<\infty$.

\begin{theorem}\label{main}
If $1\le q<p\le\infty$, then
\begin{equation}\label{mainest}
 \int_\Omega\norm{G_A}_{p\to q}\,d\prob
 \asymp_{p,q}D_1(A)+D_2(A)
\end{equation}
for every real matrix $A$, uniformly in its dimensions.
\end{theorem}

\begin{theorem}\label{mainremaining}
If $1<p\le q<2$ or $2<p\le q<\infty$, then, uniformly over finite real matrices $A$,
\begin{equation}\label{remainingest}
 \int_\Omega\norm{G_A}_{p\to q}\,d\prob
 \asymp_{p,q}D_1(A)+D_2(A)+L_{p,q}(A).
\end{equation}
\end{theorem}

For $q<p$, we use the diagonal factorization of Lemma~\ref{factor} and symmetrization of the independent rows. For $1<p\le q<2$, Lemma~\ref{rr:prelim} gives an upper bound with an additional term proportional to $\max_j\norm{(a_{ij})_i}_2[\log(2\mathrm e mn)]^3$. The proof uses the same truncation level for successive coordinate restrictions of a vector. Lemma~\ref{rr:connected} gives an estimate for matrices with sparse support; Proposition~\ref{rr:uniform} and the block decomposition in Section~\ref{sec:blocks} yield \eqref{remainingest}.

The regions of Theorem~\ref{mainremaining} are interchanged by the duality $(p,q)\mapsto(q',p')$. The remaining pairs of exponents are covered by \cite{LS} and \cite{APSS}.

\begin{theorem}\label{full}
Conjecture~1 of \cite{APSS} holds for every $1\le p,q\le\infty$. More precisely, the equivalence \eqref{conj} holds uniformly over all finite real matrices, with $L_{p,q}$ given by \eqref{lastterm}.
\end{theorem}

The second conjecture of \cite{APSS} extends to arbitrary $p$ and $q$ the characterization of bounded infinite Gaussian matrices on $\ell_2$ obtained in \cite{LVHY}. Adamczak, Prochno, Strzelecka and Strzelecki proved in \cite[Proposition 1.2]{APSS} that their first conjecture implies the second for each pair of exponents. For $a\in(0,\infty]$, let $\ell_a=\{x\in\R^{\mathbb N}:\norm{x}_a<\infty\}$, with the same formula for the (quasi-)norm as in finite dimensions, now over $\mathbb N$. An infinite matrix $C=(c_{ij})_{i,j\ge1}$ defines a bounded matrix operator $C:\ell_a\to\ell_b$ if the coordinate series $\sum_{j\ge1}c_{ij}x_j$ converge for all $i\ge1$ and $x\in\ell_a$, and the map
\[
 Cx=\left(\sum_{j\ge1}c_{ij}x_j\right)_{i\ge1}
\]
 takes values in $\ell_b$ with $\norm{C}_{a\to b}=\sup_{\norm{x}_a\le1}\norm{Cx}_b<\infty$. Otherwise set $\norm{C}_{a\to b}=\infty$. For an infinite real matrix $A=(a_{ij})_{i,j\ge1}$, take an independent family $(g_{ij})_{i,j\ge1}$ of standard Gaussian random variables and put $G_A=(a_{ij}g_{ij})_{i,j\ge1}$. Write $A_N=(a_{ij})_{1\le i,j\le N}$ for its finite principal sections.

\begin{corollary}\label{infinite}
For $1\le p,q\le\infty$, the infinite Gaussian matrix $G_A$ defines an almost surely bounded operator from $\ell_p$ to $\ell_q$ if and only if
\begin{equation}\label{infiniteest}
 \sup_{N\ge1}\bigl(D_1(A_N)+D_2(A_N)+L_{p,q}(A_N)\bigr)<\infty.
\end{equation}
Equivalently, both $A\circ A:\ell_{p/2}\to\ell_{q/2}$ and $(A\circ A)^t:\ell_{q'/2}\to\ell_{p'/2}$ are bounded, and the following additional condition holds in the indicated region:
\begin{itemize}
 \item if $p\le2\le q$, then $\int_\Omega\sup_{i,j}|a_{ij}g_{ij}|\,d\prob<\infty$;
 \item if $p\le q<2$, then $b_j\to0$ and $\sup_{j\ge1}b_j^*\sqrt{\log(j+1)}<\infty$;
 \item if $2<p\le q$, then $d_i\to0$ and $\sup_{i\ge1}d_i^*\sqrt{\log(i+1)}<\infty$.
\end{itemize}
Here $b_j$ and $d_i$ are defined by the column and row norms in \eqref{intro:colscale} and after \eqref{lastterm}, respectively, now over the infinite index sets, with $2p/(p-2)=2$ when $p=\infty$. No additional condition is required if $q<p$. In particular, Conjecture~2 of \cite{APSS} holds for all $p,q$.
\end{corollary}

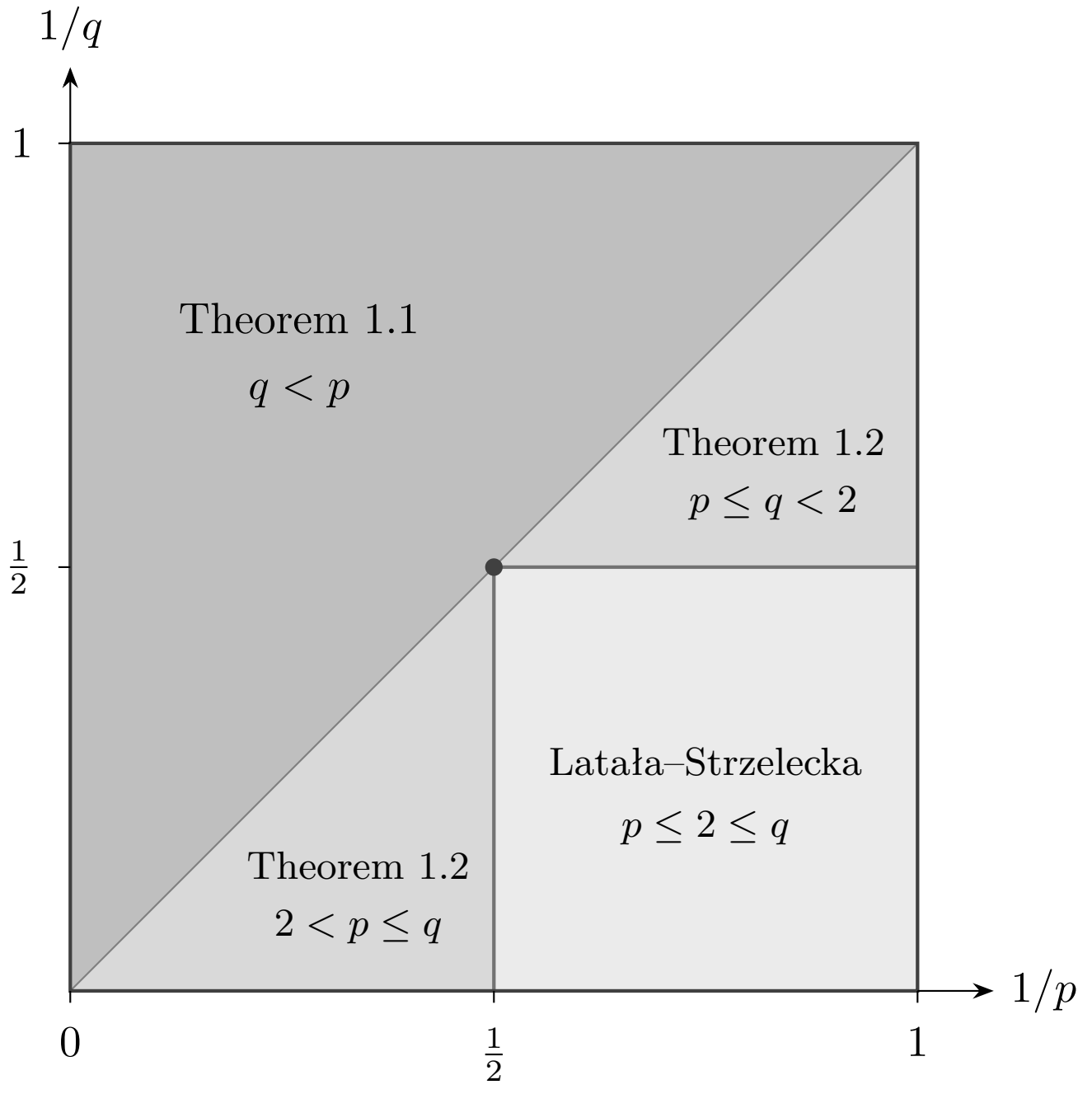
\begin{figure}[htbp]
\centering
\begin{tikzpicture}[x=6.8cm,y=6.8cm,>=Stealth,font=\small]
  \fill[black!25] (0,0)--(0,1)--(1,1)--cycle;
  \fill[black!8] (.5,0) rectangle (1,.5);
  \fill[black!15] (0,0)--(.5,0)--(.5,.5)--cycle;
  \fill[black!15] (.5,.5)--(1,.5)--(1,1)--cycle;
  \draw[black!50] (0,0)--(1,1);
  \draw[black!55,thick] (.5,0)--(.5,.5)--(1,.5);
  \draw[black!75,thick] (0,0) rectangle (1,1);
  \fill[black!75] (.5,.5) circle (2pt);
  \draw[->] (1,0)--(1.09,0) node[right] {$1/p$};
  \draw[->] (0,1)--(0,1.09) node[above] {$1/q$};
  \foreach \x/\lab in {0/0,.5/{\frac12},1/1}{\draw (\x,0)--++(0,-.014) node[below=2pt] {$\lab$};}
  \foreach \y/\lab in {.5/{\frac12},1/1}{\draw (0,\y)--++(-.014,0) node[left=2pt] {$\lab$};}
  \node[align=center] at (.27,.75) {Theorem \ref{main}\\[3pt] $q<p$};
  \node[align=center,font=\footnotesize] at (.75,.23) {Lata\l a--Strzelecka\\[3pt] $p\le2\le q$};
  \node[align=center,font=\footnotesize] at (.34,.11) {Theorem \ref{mainremaining}\\[2pt] $2<p\le q$};
  \node[align=center,font=\footnotesize] at (.83,.61) {Theorem \ref{mainremaining}\\[2pt] $p\le q<2$};
\end{tikzpicture}
\caption{The exponent square. The outer edges were treated in \cite{APSS}; the crossing region was established in \cite{LS}. The remaining regions, including the corresponding diagonal segments, are covered by Theorems~\ref{main} and \ref{mainremaining}.}
\label{fig:regions}
\end{figure}

Figure~\ref{fig:regions} distinguishes the regions treated in the present paper from those previously established.

\section{Proof of Theorem \ref{main}}\label{sec:gaussproof}
\subsection{Probabilistic notation and estimates}
All auxiliary random variables may be realized on product probability spaces. A random variable is symmetric if it has the same distribution as its negative. An independent copy of a random variable has the same distribution and is independent of that variable.

A Rademacher variable takes the values $1$ and $-1$ with probability $1/2$ each. For independent Rademacher variables $(\varepsilon_k)$ on a probability space $(\Omega_\varepsilon,\mathcal A_\varepsilon,\prob_\varepsilon)$ and any finite family of real scalars $(c_k)$, the scalar Khintchine inequality reads
\begin{equation}\label{khintchine}
 \kappa_t\Big(\sum_k c_k^2\Big)^{1/2}
 \le \left(\int\Big|\sum_k\varepsilon_k c_k\Big|^t\,d\prob_\varepsilon\right)^{1/t}
 \le \kappa_t^+\Big(\sum_k c_k^2\Big)^{1/2}\quad(t>0),
\end{equation}
where $0<\kappa_t\le \kappa_t^+<\infty$ depend only on $t$. For $t\le2$ the upper bound holds with $\kappa_t^+=1$, by monotonicity of moments and the identity $\int|\sum_k\varepsilon_k c_k|^2\,d\prob_\varepsilon=\sum_k c_k^2$. The optimal lower constant at $t=1$ is $\kappa_1=1/\sqrt{2}$ (see \cite{Haagerup}).

By the contraction principle for Rademacher sums \cite[Chapter 4]{LT}, if $T\subset\R^m$ is compact and $\varphi_i:\R\to\R$ satisfy $\varphi_i(0)=0$ and $|\varphi_i(u)-\varphi_i(v)|\le L_i|u-v|$, then
\begin{equation}\label{contraction}
\begin{split}
 \int\sup_{t\in T}\left|\sum_i\varepsilon_i\varphi_i(t_i)\right|d\prob_\varepsilon
 &\le 2\int\sup_{t\in T}\left|\sum_i\varepsilon_i L_i t_i\right|d\prob_\varepsilon.
\end{split}
\end{equation}
The suprema in \eqref{contraction} are measurable since the functions of $t\in T$ are continuous and $T$ has a countable dense subset. For a nonnegative random variable $W:\Omega\to[0,\infty)$ and $t>0$, Tonelli's theorem applied to $W^t=\int_0^W t u^{t-1}\,du$ gives
\begin{equation}\label{tailidentity}
 \int W^t\,d\prob=\int_0^\infty t u^{t-1}\prob(W>u)\,du.
\end{equation}
Integrals without subscripts refer to the probability space of the random variables in the integrand.

For $1\le q<p\le\infty$, let $\xi_{ij}:\Omega\to\R$, $1\le i\le m$, $1\le j\le n$, be independent integrable random variables with $\int_\Omega\xi_{ij}\,d\prob=0$. Set $X_A=(a_{ij}\xi_{ij}):\ell_p^n\to\ell_q^m$ and let
\begin{equation}\label{Rdef}
 s=(1/q-1/p)^{-1},\qquad R=\max\{2,s\}.
\end{equation}

\subsection{A diagonal factorization}\label{sec:factor}
In this subsection and Subsection~\ref{sec:rad}, assume that $1<q<\min\{p,2\}$ and $p<\infty$. Fix $r$ with $q<r<p$ if $p<2$, and take $r=2$ if $p\ge2$. Put
\begin{equation}\label{h}
 h=\frac{qr}{r-q},\qquad \frac1q=\frac1r+\frac1h.
\end{equation}
For $d\in\R^m$, the diagonal multiplier $D_d:\ell_r^m\to\ell_q^m$, $D_dz=(d_i z_i)_i$, satisfies
\begin{equation}\label{diagholder}
 \norm{D_d}_{r\to q}=\norm{d}_h.
\end{equation}
Indeed, H\"older with exponents $h/q$ and $r/q$ gives
\[
 \sum_i|d_i z_i|^q\le\left(\sum_i|d_i|^h\right)^{q/h}
 \left(\sum_i|z_i|^r\right)^{q/r}.
\]
If $d\ne0$, equality holds for $z_i=|d_i|^{h/r}/\norm{d}_h^{h/r}$.

The diagonal factorization below is a finite-dimensional case of Maurey's factorization theorem \cite[Th\'eor\`eme~2 and its corollary]{Maurey}. We give a direct proof.
\begin{lemma}\label{factor}
For every nonzero real matrix $V:\ell_p^n\to\ell_q^m$, there exist nonnegative numbers $d_i$, positive on the nonzero rows of $V$, such that
\begin{equation}\label{fact}
 \norm{d}_h\le C_{p,q,r}\norm{V}_{p\to q},\qquad
 \sum_{i:d_i>0}\frac{|(Vx)_i|^r}{d_i^r}\le1\quad(x\in B_p^n).
\end{equation}
\end{lemma}
\begin{proof}
Discard the zero rows, and let $K$ be the convex hull in $\R^m$ of
\[
 \{(|(Vx)_i|^r)_{i=1}^m:x\in B_p^n\}.
\]
This is a compact set. Define $F=\max_{z\in K}(\sum_i z_i^{q/r})^{1/q}$. We have
\begin{equation}\label{Fbound}
 F\le C_{p,q,r}\norm{V}_{p\to q}.
\end{equation}
Each $z\in K$ can be written as $z_i=\sum_k\theta_k|(Vx_k)_i|^r$, with $\theta_k\ge0$, $\sum_k\theta_k=1$ and $x_k\in B_p^n$.

If $p\ge2$, let $g_k$ be independent standard Gaussians and set $Z=\sum_k\theta_k^{1/2}g_kx_k$. Since $r=2$, Minkowski's inequality gives
\begin{align*}
 \int\norm{Z}_p^p\,d\prob
 &=\gamma_p^p\sum_j\left(\sum_k\theta_k|(x_k)_j|^2\right)^{p/2}\\
 &\le\gamma_p^p\left(\sum_k\theta_k\norm{x_k}_p^2\right)^{p/2}
 \le\gamma_p^p.
\end{align*}
Consequently,
\begin{align*}
 \gamma_q^q\sum_i z_i^{q/2}
 &=\int\norm{VZ}_q^q\,d\prob\\
 &\le\norm{V}_{p\to q}^q\left(\int\norm{Z}_p^p\,d\prob\right)^{q/p}
 \le\gamma_p^q\norm{V}_{p\to q}^q.
\end{align*}

If $p<2$, the exponent $r$ lies strictly between $q$ and $p$. Take independent uniform random variables $U_k$ on $(0,1)$ and independent Rademacher variables $\varepsilon_k$, independent of the $U_k$. The symmetric random variables $\xi_k=\varepsilon_k U_k^{-1/r}$ satisfy $\prob(|\xi_k|>t)=\min\{1,t^{-r}\}$ for $t>0$. Their absolute moments of order $q<r$ are finite. For any finite scalar family $(c_k)$, put $b=(\sum_k|c_k|^r)^{1/r}$. If $b>0$, independence gives
\[
 \prob\left(\max_k|c_k\xi_k|>b\right)
 \ge 1-\exp\left(-\sum_k\frac{|c_k|^r}{b^r}\right)=1-e^{-1}.
\]
Condition on the absolute values of the $\xi_k$; the signs are independent Rademacher variables. By convexity, averaging over all the signs except one corresponding to a largest term shows that the conditional moment of order $q$ is at least the $q$-th power of that modulus. Hence
\begin{equation}\label{heavylower}
 \int\left|\sum_k c_k\xi_k\right|^q\,d\prob
 \ge (1-e^{-1})\left(\sum_k|c_k|^r\right)^{q/r}.
\end{equation}
Put $Y_k=\theta_k^{1/r}\xi_k x_k$ and $Z=\sum_kY_k$. For each $t>0$,
\[
 \sum_k\prob(\norm{Y_k}_p>t)\le t^{-r}\sum_k\theta_k\norm{x_k}_p^r\le t^{-r}.
\]
The truncated vectors $Y_k^{(t)}=Y_k\mathbf1_{\{\norm{Y_k}_p\le t\}}$ remain independent and symmetric. Since each truncated vector is symmetric and the vectors are independent, their joint distribution is unchanged by independent signs. Applying \eqref{khintchine} in each coordinate, with $\kappa_p^+=1$, and then \eqref{tailidentity} to the truncated norms gives
\begin{align*}
 \int\norm{\sum_kY_k^{(t)}}_p^p\,d\prob
 &\le\sum_j\int\left(\sum_k|(Y_k^{(t)})_j|^2\right)^{p/2}d\prob\\
 &\le\sum_k\int\norm{Y_k^{(t)}}_p^p\,d\prob\\
 &\le\sum_k\int_0^t p u^{p-1}\prob(\norm{Y_k}_p>u)\,du\\
 &\le\frac{p}{p-r}t^{p-r}.
\end{align*}
On the event $\bigcap_k\{\norm{Y_k}_p\le t\}$, one has $Z=\sum_kY_k^{(t)}$. The union bound and Markov's inequality give
\begin{align*}
 \prob(\norm{Z}_p>t)
 &\le\sum_k\prob(\norm{Y_k}_p>t)
       +\prob\left(\norm{\sum_kY_k^{(t)}}_p>t\right)\\
 &\le t^{-r}+t^{-p}\int\norm{\sum_kY_k^{(t)}}_p^p\,d\prob\\
 &\le\left(1+\frac{p}{p-r}\right)t^{-r}.
\end{align*}
Let $C=1+p/(p-r)$. Using \eqref{tailidentity} and the trivial bound $\prob(\norm{Z}_p>t)\le1$, we obtain
\[
 \int\norm{Z}_p^q\,d\prob\le \int_0^1 qt^{q-1}\,dt+C\int_1^\infty qt^{q-1-r}\,dt
 =1+\frac{Cq}{r-q}.
\]
By \eqref{heavylower}, applied to each coordinate of $VZ$,
\[
 (1-e^{-1})\sum_i z_i^{q/r}\le\int\norm{VZ}_q^q\,d\prob
 \le C_{p,q,r}\norm{V}_{p\to q}^q.
\]

Let $z^*\in K$ maximize $\sum_i z_i^{q/r}$. Since each row of $V$ is nonzero, a convex combination of elements of $K$ is strictly positive. Let $w\in K$ have strictly positive coordinates, and suppose that $J_0=\{i:z_i^*=0\}$ is nonempty. For $0<t\le1/2$, the vector $(1-t)z^*+tw$ belongs to $K$ and
\begin{align*}
 &\sum_i((1-t)z_i^*+tw_i)^{q/r}-\sum_i(z_i^*)^{q/r}\\
 &\quad\ge t^{q/r}\sum_{i\in J_0}w_i^{q/r}
 +\big((1-t)^{q/r}-1\big)\sum_{i\notin J_0}(z_i^*)^{q/r}\\
 &\quad\ge t^{q/r}\sum_{i\in J_0}w_i^{q/r}-C_{q,r}tF^q.
\end{align*}
Since $q/r<1$, this is positive for sufficiently small $t$, contradicting maximality. Thus $z_i^*>0$ for every $i$. For $z\in K$ and $0\le t\le1$, $(1-t)z^*+tz\in K$. Hence
\begin{align*}
 0&\ge\left.\frac{d}{dt}\sum_i((1-t)z_i^*+tz_i)^{q/r}\right|_{t=0}\\
 &=\frac qr\sum_i(z_i^*)^{q/r-1}(z_i-z_i^*).
\end{align*}
It follows that
\begin{equation}\label{factorstationary}
 \sum_i(z_i^*)^{q/r-1}z_i\le\sum_i(z_i^*)^{q/r}=F^q.
\end{equation}
Define $d_i^r=F^q(z_i^*)^{1-q/r}$. By \eqref{factorstationary}, for $x\in B_p^n$,
\[
 \sum_i\frac{|(Vx)_i|^r}{d_i^r}
 =F^{-q}\sum_i(z_i^*)^{q/r-1}|(Vx)_i|^r\le1.
\]
Since $h/r=q/(r-q)$, we have
\[
 (1-q/r)h/r=q/r,\qquad qh/r+q=h.
\]
Consequently,
\begin{align*}
 \sum_i d_i^h&=F^{qh/r}\sum_i(z_i^*)^{(1-q/r)h/r}\\
 &=F^{qh/r}\sum_i(z_i^*)^{q/r}=F^{qh/r+q}=F^h.
\end{align*}
Thus $\norm{d}_h=F\le C_{p,q,r}\norm{V}_{p\to q}$ by \eqref{Fbound}. Set $d_i=0$ on the discarded rows.
\end{proof}

\subsection{A Rademacher estimate}\label{sec:rad}
Let $V=(v_{ij}):\ell_p^n\to\ell_q^m$ be a deterministic real matrix and $v_i=(v_{ij})_{j=1}^n\in\R^n$ its rows. Write $M=\norm{V}_{p\to q}$ and let $(\varepsilon_i)$ be independent Rademacher variables on a probability space $(\Omega_\varepsilon,\mathcal A_\varepsilon,\prob_\varepsilon)$. Define
\[
 S(V)=\int_{\Omega_\varepsilon}\sup_{x\in B_p^n}
 \left|\sum_i\varepsilon_i |(Vx)_i|^q\right|\,d\prob_\varepsilon.
\]
\begin{lemma}\label{rad}
With $r$ as in Subsection~\ref{sec:factor}, set $\alpha=r(q-1)/(r-1)$ and $\beta=(r-q)/(r-1)$. Then
\begin{equation}\label{radbound}
 S(V)\le C_{p,q,r}M^\alpha D_2(V)^\beta.
\end{equation}
\end{lemma}
\begin{proof}
The case $V=0$ is immediate. Otherwise take the factors $d_i$ from Lemma \ref{factor}, omitting zero rows. For $\lambda>0$, let
\[
 L_i=d_i^{r/(r-q)}\lambda^{-1/(r-q)},\qquad
 \varphi_i(t)=\min\{|t|^q,L_i^q\}\quad(t\in\R).
\]
For $|t|\le L_i$, the difference $|t|^q-\varphi_i(t)$ is zero. For $|t|>L_i$, since $q-r<0$,
\[
 0\le |t|^q-L_i^q\le |t|^r|t|^{q-r}
 \le |t|^r L_i^{q-r}=\lambda |t|^r/d_i^r.
\]
By \eqref{fact},
\begin{equation}\label{err}
 \sup_{x\in B_p^n}\sum_i\big(|(Vx)_i|^q-\varphi_i((Vx)_i)\big)\le\lambda.
\end{equation}
The derivative of $t\mapsto |t|^q$ has absolute value at most $qL_i^{q-1}$ on $[-L_i,L_i]$. Outside this interval $\varphi_i$ is constant. Hence $\varphi_i(0)=0$ and its Lipschitz constant is at most $qL_i^{q-1}$. Apply \eqref{contraction} with $T=V(B_p^n)$. Since $\sup_{x\in B_p^n}|\langle b,x\rangle|=\norm{b}_{p'}$, we have
\begin{align*}
 &\int\sup_{x\in B_p^n}\left|\sum_i\varepsilon_i\varphi_i((Vx)_i)\right|d\prob_\varepsilon\\
 &\qquad\le 2q\int\norm{\sum_i\varepsilon_i L_i^{q-1}v_i}_{p'}d\prob_\varepsilon.
\end{align*}
For any nonnegative vector $w$, the scalar Khintchine inequality and Jensen's inequality give
\begin{align*}
 \int\norm{\sum_i\varepsilon_iw_iv_i}_{p'}d\prob_\varepsilon
 &\le \kappa_{p'}^+\left[\sum_j\left(\sum_i w_i^2v_{ij}^2\right)^{p'/2}\right]^{1/p'}\\
 &\le \kappa_{p'}^+D_2(V)\norm{w}_{q'}.
\end{align*}
Since $(q-1)q'=q$, \eqref{h} gives $r(q-1)q'/(r-q)=h$. Thus
\begin{align*}
 \norm{(L_i^{q-1})_i}_{q'}
 &=\left(\sum_i d_i^{r(q-1)q'/(r-q)}
 \lambda^{-(q-1)q'/(r-q)}\right)^{1/q'}\\
 &=\lambda^{-(q-1)/(r-q)}\left(\sum_i d_i^h\right)^{1/q'}\\
 &=\lambda^{-(q-1)/(r-q)}\norm{d}_h^{r(q-1)/(r-q)}.
\end{align*}
Using \eqref{fact} and \eqref{err}, we obtain
\begin{equation}\label{lambda}
 S(V)\le\lambda+C_{p,q,r}D_2(V)M^{r(q-1)/(r-q)}\lambda^{-(q-1)/(r-q)}.
\end{equation}
Write $b=(q-1)/(r-q)>0$. Then $\alpha=rb/(1+b)$ and $\beta=1/(1+b)$. Choose $\lambda=M^\alpha D_2(V)^\beta>0$. The second term on the right of \eqref{lambda}, apart from the constant, is
\[
 D_2(V)M^{rb}(M^\alpha D_2(V)^\beta)^{-b}
 =M^{rb-\alpha b}D_2(V)^{1-\beta b}
 =M^\alpha D_2(V)^\beta.
\]
Finally,
\[
 \alpha+\beta=\frac{r(q-1)+r-q}{r-1}=q,
 \qquad \beta>0.
\]
\end{proof}

\subsection{Random perturbations of positive matrices}
\begin{lemma}\label{positive}
Let $1<a<\infty$ and $0<b<a$. Let $B$ be a nonnegative real matrix of size $m\times n$ and let $Z_{ij}:\Omega\to[0,\infty)$, $1\le i\le m$, $1\le j\le n$, be random variables. Write $Z=(Z_{ij})$ and view $B\circ Z$ as a random linear map $\ell_a^n\to\ell_b^m$. Set $v=\max\{a',ab/(a-b)\}$. Then
\begin{equation}\label{positiveest}
 \int\norm{B\circ Z}_{a\to b}\,d\prob
 \le\norm{B}_{a\to b}\max_{i,j}\left(\int Z_{ij}^v\,d\prob\right)^{1/v}.
\end{equation}
\end{lemma}
\begin{proof}
It suffices to consider matrices with strictly positive entries: replace $B$ by $B+\epsilon\mathbf 1$, use monotonicity of the norms on positive matrices and let $\epsilon\downarrow0$ by finite-dimensional continuity. Normalize $\norm{B}_{a\to b}=1$. There is a nonnegative vector $x$ of $\ell_a$-norm one at which this norm is attained. If $x_j=0$, then for $t>0$ the vector $x(t)=(x+te_j)/(1+t^a)^{1/a}$ has norm one. Since $Bx$ has strictly positive coordinates,
\[
 \lim_{t\downarrow0}\frac{\norm{Bx(t)}_b^b-\norm{Bx}_b^b}{t}
 =b\sum_i b_{ij}(Bx)_i^{b-1}>0,
\]
where $a>1$ ensures that the derivative of $(1+t^a)^{-b/a}$ at zero is zero. This contradicts maximality; hence $x_j>0$ for every $j$.

Write $y=Bx$, so that $\sum_i y_i^b=1$. The Lagrange multiplier equation for maximizing $\sum_i(Bx)_i^b$ under $\sum_j x_j^a=1$ is
\[
 B^t(y_i^{b-1})_i=(x_j^{a-1})_j.
\]
Indeed, writing the multiplier as $\Lambda$, the coordinate equations are
\[
 b\sum_i b_{ij}y_i^{b-1}=\Lambda a x_j^{a-1}.
\]
Multiplication by $x_j$ and summation over $j$ give
\[
 \Lambda a=b\sum_i y_i^{b-1}\sum_j b_{ij}x_j
 =b\sum_i y_i^b=b.
\] Put
\[
 \mu_j=x_j^a,\qquad \nu_i=y_i^b,\qquad P_{ij}=\frac{b_{ij}x_j}{y_i}.
\]
The vectors $\mu,\nu$ are probability vectors, $\sum_jP_{ij}=1$, and the multiplier equation gives
\begin{equation}\label{transport}
 \sum_i\nu_iP_{ij}=x_j\sum_i b_{ij}y_i^{b-1}=\mu_j.
\end{equation}
Every vector $z\in\R^n$ has the form $z_j=x_jf_j$. By H\"older,
\[
 \left|\sum_j P_{ij}Z_{ij}f_j\right|
 \le U_i W_i,
\]
where
\[
 U_i=\left(\sum_jP_{ij}|f_j|^a\right)^{1/a},\quad
 W_i=\left(\sum_jP_{ij}Z_{ij}^{a'}\right)^{1/a'}.
\]
Let $u=ab/(a-b)$. Since $1/b=1/a+1/u$, H\"older with the exponents $a/b,u/b>1$ gives
\begin{align*}
 \norm{(B\circ Z)z}_b
 &\le\left(\sum_i\nu_i U_i^a\right)^{1/a}
       \left(\sum_i\nu_i W_i^u\right)^{1/u}\\
 &\le\norm{z}_a\left(\sum_{i,j}\nu_i P_{ij}Z_{ij}^v\right)^{1/v}.
\end{align*}
In fact, \eqref{transport} gives
\begin{align*}
 \sum_i\nu_i U_i^a&=\sum_j|f_j|^a\sum_i\nu_iP_{ij}\\
 &=\sum_j x_j^a|f_j|^a=\norm{z}_a^a.
\end{align*}
Since $v\ge a'$ and $\sum_jP_{ij}=1$,
\[
 W_i\le\left(\sum_jP_{ij}Z_{ij}^v\right)^{1/v}.
\]
Since also $v\ge u$ and $\sum_i\nu_i=1$,
\begin{align*}
 \left(\sum_i\nu_i W_i^u\right)^{1/u}
 &\le\left(\sum_i\nu_i W_i^v\right)^{1/v}\\
 &\le\left(\sum_{i,j}\nu_i P_{ij} Z_{ij}^v\right)^{1/v}.
\end{align*}
Taking the supremum over $z\in B_a^n$ and applying Jensen's inequality to the concave function $t^{1/v}$ gives
\[
 \int\norm{B\circ Z}_{a\to b}\,d\prob
 \le\left(\sum_{i,j}\nu_iP_{ij}\int Z_{ij}^v\,d\prob\right)^{1/v}
 \le\max_{i,j}\norm{Z_{ij}}_{L_v}.
\]
Restoring the normalization proves \eqref{positiveest}.
\end{proof}

\subsection{Direct upper estimates}\label{direct}
\begin{lemma}\label{directupper}
Let $1\le q<p\le\infty$. Set
\begin{equation}\label{directmoment}
 R_* = \begin{cases}
 \max\{R,2q/(2-q)\},&1<q<2<p<\infty,\\
 R,&\text{otherwise},
 \end{cases}
 \qquad K_* = \max_{i,j}\norm{\xi_{ij}}_{L_{R_*}}.
\end{equation}
For independent centered real entries with $K_*<\infty$,
\begin{equation}\label{auxupper}
 \int\norm{X_A}_{p\to q}\,d\prob\le C_{p,q}K_*(D_1(A)+D_2(A)).
\end{equation}
\end{lemma}
\begin{proof}
First suppose $1<q<\min\{p,2\}$ and $p<\infty$. For $x\in B_p^n$, independence and centering give
\begin{align}
 \int\left|\sum_j a_{ij}\xi_{ij}x_j\right|^q d\prob
 &\le\left(\int\left|\sum_j a_{ij}\xi_{ij}x_j\right|^2d\prob\right)^{q/2}\notag\\
 &\le K_*^q\left(\sum_j a_{ij}^2x_j^2\right)^{q/2}.\label{ngvariance}
\end{align}
Apply Lemma \ref{positive} with $B=(A\circ A)^t$, $Z=(\xi_{ij}^2)^t$, $a=q'/2$ and $b=p'/2$. In this case
\begin{equation}\label{momentalgebra}
 2a'=\frac{2q}{2-q},\qquad
 \frac{2ab}{a-b}=\frac{pq}{p-q}=s.
\end{equation}
The exponent $2v$ in Lemma \ref{positive} is therefore $\max\{2q/(2-q),s\}$. If $p\le2$, then $1/s=1/q-1/p\le1/q-1/2$, so $s\ge2q/(2-q)>2$ and $2v=R$. If $p>2$, then $2v=R_*$ by \eqref{directmoment}. In both cases $2v=R_*$. Consequently,
\begin{align*}
 \int D_2(X_A)^2\,d\prob
 &=\int\norm{(A\circ A)^t\circ(\xi_{ij}^2)^t}_{q'/2\to p'/2}\,d\prob\\
 &\le D_2(A)^2\max_{i,j}\left(\int|\xi_{ij}|^{R_*}\,d\prob\right)^{2/R_*}\\
 &=K_*^2 D_2(A)^2.
\end{align*}
Since $q<2$,
\begin{equation}\label{ngd2}
 \norm{D_2(X_A)}_{L_q}\le\norm{D_2(X_A)}_{L_2}\le K_* D_2(A).
\end{equation}
Let $U=(\int\norm{X_A}_{p\to q}^q\,d\prob)^{1/q}$ and define the random functions $f_i:\Omega\times B_p^n\to[0,\infty)$ by $f_i(\omega,x)=|(X_A(\omega)x)_i|^q$. We suppress $\omega$ in the notation. The functions $f_i$ are independent since they depend on distinct rows. The fluctuation of their sum around its mean is controlled by symmetrization: introduce a second independent realization of the entire random matrix and denote the corresponding functions by $f_i'$. For each fixed realization of the original matrix, $\int f_i'(x)\,d\prob'=\int f_i(x)\,d\prob$. Thus Jensen's inequality on the second probability space bounds the absolute centered sum by the integral of the absolute difference. Taking the supremum and using \eqref{ngvariance}, we obtain
\begin{align*}
 U^q&\le K_*^qD_1(A)^q+\int\sup_{x\in B_p^n}\left|\sum_i(f_i(x)-\int f_i(x)\,d\prob)\right|d\prob\\
 &\le K_*^qD_1(A)^q+\iint\sup_{x\in B_p^n}\left|\sum_i(f_i(x)-f_i'(x))\right|d\prob' d\prob\\
 &\le K_*^qD_1(A)^q+2\int S(X_A)\,d\prob.
\end{align*}
Indeed, the independence of the pairs $(f_i,f_i')$ and invariance under interchanging the two entries of each pair imply
\begin{align*}
 &\iint\sup_x\left|\sum_i(f_i(x)-f_i'(x))\right|\,d\prob' d\prob\\
 &=\iiint\sup_x\left|\sum_i\varepsilon_i(f_i(x)-f_i'(x))\right|
 \,d\prob_\varepsilon d\prob' d\prob\\
 &\le\iiint\left(\sup_x\left|\sum_i\varepsilon_i f_i(x)\right|
 +\sup_x\left|\sum_i\varepsilon_i f_i'(x)\right|\right)
 \,d\prob_\varepsilon d\prob' d\prob\\
 &=2\int S(X_A)\,d\prob,
\end{align*}
where every supremum is over $B_p^n$. All suprema may be taken over countable dense subsets of the unit ball; the integrals are finite by the moment assumption and finite dimensionality. By Lemma \ref{rad}, H\"older with exponents $q/\alpha,q/\beta$ and \eqref{ngd2},
\begin{align}
 U^q&\le K_*^qD_1(A)^q+C_{p,q,r}\int\norm{X_A}_{p\to q}^\alpha D_2(X_A)^\beta\,d\prob\notag\\
 &\le K_*^qD_1(A)^q+C_{p,q,r} U^\alpha (K_*D_2(A))^\beta.\label{absorb}
\end{align}
If $U^q\le 2K_*^qD_1(A)^q$, then $U\le 2^{1/q}K_*D_1(A)$. Otherwise \eqref{absorb} gives $U^\beta\le 2C_{p,q,r}(K_*D_2(A))^\beta$, since $\alpha+\beta=q$ and $\beta>0$. This proves \eqref{auxupper} for $1<q<\min\{p,2\}$ and $p<\infty$, with $r=(p+q)/2$ for $p<2$ and $r=2$ for $p\ge2$. If $2\le q<p<\infty$, then $1<p'<q'\le2$. For any real matrix $V$, duality gives
\begin{equation}\label{duality}
 \norm{V}_{p\to q}=\sup_{x\in B_p^n,\,y\in B_{q'}^m}|\langle Vx,y\rangle|=\norm{V^t}_{q'\to p'}.
\end{equation}
For the transposed problem the difference of reciprocal exponents is $1/p'-1/q'=1/q-1/p$, and $D_1,D_2$ are interchanged by their definitions. Both exponent pairs lie outside the first case in \eqref{directmoment}; hence $R_*=R$ is unchanged. Applying \eqref{auxupper} with exponents $(q',p')$ to $X_A^t$ proves this case.

For $q=1<p\le\infty$, Cauchy--Schwarz and the independence and centering of the entries give $\int_\Omega |(X_Ax)_i|\,d\prob\le K_*(\sum_j a_{ij}^2x_j^2)^{1/2}$. Summing over $i$ bounds the sum of the first absolute moments by $K_*D_1(A)$ for every $x\in B_p^n$. Symmetrization of the independent rows, contraction with the function $t\mapsto |t|$, and the scalar Khintchine inequality give
\begin{align}
 \int\norm{X_A}_{p\to1}\,d\prob
 &\le K_* D_1(A)+4\iint\norm{\sum_i\varepsilon_i (a_{ij}\xi_{ij})_{j=1}^n}_{p'}\,d\prob_\varepsilon d\prob\notag\\
 &\le K_* D_1(A)+4\kappa_{p'}^+\int D_2(X_A)\,d\prob.\label{qone}
\end{align}
Here $p'<\infty$ and $D_2(V)=(\sum_j(\sum_i v_{ij}^2)^{p'/2})^{1/p'}$ because $q'=\infty$. In this case $R=\max\{2,p'\}$. If $p'\le2$, Jensen's inequality implies
\[
 \int\left(\sum_i a_{ij}^2\xi_{ij}^2\right)^{p'/2}d\prob
 \le K_*^{p'}\left(\sum_i a_{ij}^2\right)^{p'/2}.
\]
If $p'>2$, then $R=p'$ and Minkowski's inequality gives
\begin{align*}
 \left(\int\left(\sum_i a_{ij}^2\xi_{ij}^2\right)^{p'/2}d\prob\right)^{2/p'}
 &\le\sum_i a_{ij}^2\left(\int|\xi_{ij}|^{p'}\,d\prob\right)^{2/p'}\\
 &\le K_*^2\sum_i a_{ij}^2.
\end{align*}
In either case, Jensen's inequality gives
\begin{align*}
 \int D_2(X_A)\,d\prob
 &\le\left(\sum_j\int\left(\sum_i a_{ij}^2\xi_{ij}^2\right)^{p'/2}
 d\prob\right)^{1/p'}\\
 &\le K_*\left(\sum_j\left(\sum_i a_{ij}^2\right)^{p'/2}\right)^{1/p'}
 =K_* D_2(A).
\end{align*} Substitution in \eqref{qone} gives \eqref{auxupper} for $q=1$; \eqref{duality} gives the case $p=\infty$.

\end{proof}

\subsection{The Gaussian estimate}
For Gaussian entries, $K_* = \gamma_{R_*}<\infty$ for every pair $1\le q<p\le\infty$. Lemma \ref{directupper} gives
\begin{equation}\label{gaussdirect}
 \int_\Omega\norm{G_A}_{p\to q}\,d\prob
 \le C_{p,q}\gamma_{R_*}(D_1(A)+D_2(A)).
\end{equation}
The reverse estimate is \eqref{intro:basiclower}.

\section{Proof of Theorem \ref{mainremaining}}\label{sec:remaining}
Throughout this section assume $1<p\le q<2$, and set
\begin{equation}\label{rr:notation}
 t=p',\qquad v=q',\qquad s=\frac{2q}{2-q},\qquad
 b_j=\norm{(a_{ij})_{i=1}^m}_s.
\end{equation}
Write $b_*=\max_j b_j$, $B(A)=\max_j b_j^*\sqrt{\log(j+1)}$ and $S(A)=D_1(A)+D_2(A)+B(A)$. Thus $B(A)=L_{p,q}(A)$ in this region. Write $U\lesssim_{p,q}V$ if $U\le C_{p,q}V$ for a constant independent of the dimensions and coefficient matrix; other subscripts indicate additional permitted dependencies.

By \eqref{d1},
\begin{equation}\label{rr:d1}
 D_1(A)=\max_j\norm{(a_{ij})_{i=1}^m}_q.
\end{equation}
Indeed, for $x\in B_p^n$, since $q/2<1$,
\[
 \sum_i\Big(\sum_j a_{ij}^2x_j^2\Big)^{q/2}
 \le\sum_{i,j}|a_{ij}x_j|^q
 \le\max_j\sum_i|a_{ij}|^q\,\norm{x}_q^q
 \le\max_j\sum_i|a_{ij}|^q.
\]
Equality in the supremum is obtained by testing a coordinate vector. Define $u:\R^n\to[0,\infty)^m$ and $\widetilde u:\R^m\to[0,\infty)^n$ by
\[
 u_i(x)=\Big(\sum_{j=1}^n a_{ij}^2x_j^2\Big)^{1/2},\qquad
 \widetilde u_j(y)=\Big(\sum_{i=1}^m a_{ij}^2y_i^2\Big)^{1/2}.
\]
For fixed deterministic vectors $x\in\R^n$, $y\in\R^m$, Jensen's inequality and the Gaussian moment formula \eqref{fixedx} give
\begin{equation}\label{rr:fixed}
 \int_\Omega\norm{G_Ax}_q\,d\prob\le\gamma_q\norm{u(x)}_q,\qquad
 \int_\Omega\norm{G_A^ty}_t\,d\prob\le\gamma_t\norm{\widetilde u(y)}_t.
\end{equation}
The respective suprema of the deterministic norms on the unit balls are $D_1(A)$ and $D_2(A)$. The lower bound in \eqref{remainingest} follows from \eqref{intro:basiclower} and the ordered column estimate following \eqref{intro:colscale}. 

By the Gaussian concentration inequality \cite{LT}, if $F:\R^N\to\R$ is $\lambda$-Lipschitz for the Euclidean norm and $g:\Omega\to\R^N$ is a standard Gaussian vector, then
\begin{equation}\label{rr:conc}
 \prob\big(F(g)-\int_\Omega F(g)\,d\prob>u\big)\le \exp(-u^2/(2\lambda^2))\qquad(u>0).
\end{equation}
For the matrix norm and all suprema of its bilinear evaluations over subsets of $B_p^n\times B_v^m$, the Lipschitz constant is at most $b_*$. Indeed, $1/2=1/s+1/v$ and $p<2$ imply
\[
 \sum_{i,j}a_{ij}^2x_j^2y_i^2
 \le\sum_j x_j^2 b_j^2\norm{y}_v^2
 \le b_*^2\norm{x}_p^2\norm{y}_v^2.
\]

\subsection[ A common truncation level]{A common truncation level}
For a real number $u$, put $u_+=\max\{u,0\}$. For $N,L\in\mathbb N$, define $T_L:\R^N\to[0,\infty)$ by $T_L(z)=\sum_{j=1}^{L}|z|_j^*$, where $(|z|_j^*)$ is the nonincreasing rearrangement of $(|z_j|)$, extended by zeros. Fix $h=2+\max(q,t)$. For $u\in[0,\infty)^N$, define $H_L:[0,\infty)^N\to[0,\infty)$ by
\[
 H_L(u)=\sum_{j\le L}u_j^*+
  L^{1-1/h}\left(\sum_{j>L}(u_j^*)^h\right)^{1/h}.
\]
For $1<c<h$ and nonnegative finite vectors $u,z$, extended by zero after their last coordinates,
\begin{equation}\label{rr:hardy}
 \sum_{l\ge0} z_{2^l}^* H_{2^l}(u)
       \lesssim_{c,h}\norm{z}_{c'}\norm{u}_c .
\end{equation}
Indeed, put $U_k=2^{k/c}u_{2^k}^*$. The intervals $[2^k,2^{k+1})$ contain $2^k$ integers, on each of which $u_j^*\le u_{2^k}^*$. Hence
\begin{align*}
 \sum_{j\le 2^l}u_j^*&\le\sum_{k\le l}2^k u_{2^k}^*,\\
 \left(\sum_{j>2^l}(u_j^*)^h\right)^{1/h}
 &\le\left(\sum_{k\ge l}2^k(u_{2^k}^*)^h\right)^{1/h}
 \le\sum_{k\ge l}2^{k/h}u_{2^k}^*.
\end{align*}
Multiplying the first estimate by $2^{-l/c'}$ and the second by $2^{l(1-1/h-1/c')}$ bounds $2^{-l/c'}H_{2^l}(u)$ by
\[
 \sum_{k\le l}2^{-(l-k)/c'}U_k+
 \sum_{k\ge l}2^{-(k-l)(1/c-1/h)}U_k .
\]
Both kernels $\big(2^{-j/c'}\big)_{j\ge0}$ and $\big(2^{-j(1/c-1/h)}\big)_{j\ge0}$ belong to $\ell_1$. Young's inequality for convolution therefore bounds the $\ell_c$-norm of $(2^{-l/c'}H_{2^l}(u))_{l\ge0}$ by $C_{c,h}\norm{U}_c$. Moreover,
\[
 \sum_{k\ge0}2^k(u_{2^k}^*)^c\le (u_1^*)^c+2\sum_{k\ge1}\sum_{j=2^{k-1}+1}^{2^k}(u_j^*)^c\le 2\norm{u}_c^c.
\]
The same inequality with $c'$ and $z$ bounds the $\ell_{c'}$-norm of $(2^{l/c'}z_{2^l}^*)_{l\ge0}$. H\"older's inequality for the product of these two sequences proves \eqref{rr:hardy}.

\begin{lemma}\label{rr:prelim}
Let $1<p\le q<2$ and $C_2=\max_j\norm{(a_{ij})_i}_2$. Then
\[
 \int_\Omega\norm{G_A}_{p\to q}\,d\prob\lesssim_{p,q}
 D_1(A)+D_2(A)+C_2\,[\log(2\mathrm e mn)]^3 .
\]
\end{lemma}
\begin{proof}
If $A=0$, there is nothing to prove. By homogeneity suppose $M=\max_{ij}|a_{ij}|=1$, so that $C_2\ge1$, and put $D=2\mathrm e mn$. Let $(a,c,N,M_1,Z)$ be either $(p,q,n,m,G_A)$ or $(v,t,m,n,G_A^t)$, so that $Z:\ell_a^N\to\ell_{c}^{M_1}$. Define $\sigma:\R^N\to[0,\infty)^{M_1}$ by $\sigma_i(w)=(\int_\Omega|Z_i(w)|^2\,d\prob)^{1/2}$, $1\le i\le M_1$. For $\lambda\ge0$, define $\phi_\lambda:\R\to[0,\infty)$, $\Psi_\lambda:\Omega\times\R^N\to[0,\infty)$ and $\psi_\lambda:\R^N\to[0,\infty)$ by
\[
 \phi_\lambda(z)=(|z|-\lambda)_+,\qquad
 \Psi_\lambda(w)=\sum_i\phi_\lambda(Z_i(w)),\qquad
 \psi_\lambda(w)=\int_\Omega\Psi_\lambda(w)\,d\prob.
\]
With the probability variable suppressed in $\Psi_\lambda$, for every $w\in\R^N$,
$T_L(Z(w))\le L\lambda+\Psi_\lambda(w)$.

\smallskip
\noindent\emph{Choice of the threshold.}
Fix $w\in B_a^N$ and a positive integer $L$. Set $u=\sigma(w)$ and choose
\[
 \lambda_0=\gamma_{2h}L^{-1/h}
           \left(\sum_{i>L}(u_i^*)^h\right)^{1/h}.
\]
Round $\lambda_0$ upwards to a multiple $\lambda$ of $D^{-4}$, leaving zero at zero. Every $z\in\R^N$ with $|z_j|\le |w_j|$ for $1\le j\le N$ has coordinate variances bounded by those of $w$. Separating the $L$ largest coordinates of $u$ gives
\begin{equation}\label{rr:activity}
 \begin{split}
 L\lambda+\psi_\lambda(z)&\le C_h H_L(u)+L D^{-4},\\
 \sum_i\prob\{|Z_i(z)|>\lambda\}^{1/2}&\le 2L .
 \end{split}
\end{equation}
To prove the first inequality when $\lambda>0$, let $I$ consist of the indices of the $L$ largest coordinates of $u$ (or all indices if fewer are available). On $I$, $\int\phi_\lambda(Z_i(z))\,d\prob\le\gamma_1u_i$. On its complement,
\begin{align*}
 \sum_{i\notin I}\int\phi_\lambda(Z_i(z))\,d\prob
 &\le\lambda^{1-h}\gamma_h^h\sum_{i\notin I}u_i^h\\
 &\le \frac{\gamma_h^h}{\gamma_{2h}^{h-1}} L^{1-1/h}
       \left(\sum_{i>L}(u_i^*)^h\right)^{1/h},
\end{align*}
where a zero tail is treated separately. Also $L\lambda\le L\lambda_0+LD^{-4}$. These estimates give the asserted first bound. For the second, Markov's inequality at exponent $2h$ bounds the sum outside the deterministic top $L$ by
\[
 \gamma_{2h}^h\lambda^{-h}\sum_{i>L}(u_i^*)^h\le L.
\] If $\lambda=0$, at most $L$ coordinates have nonzero variance and both conclusions follow directly. These estimates hold at the same threshold for all coordinate restrictions and inward roundings of $w$.

\smallskip
\noindent\emph{Increments on a finite grid.}
In either input space round coordinates towards zero to multiples of $D^{-4}$. Consider any two consecutive rounded prefixes $w_0,w_1=w_0+e$ ending with an input shell of ranks $K,\ldots,2K-1$, with $K$ a power of two. The possible pairs are contained in the collection having at most $2K$ coordinates in the union of their supports, entries in the grid on $[-1,1]$, and disjoint supports for $w_0,e$. Each coordinate of the grid has at most $2D^4+1$ possible values. If $N$ is the input dimension, the number of pairs is bounded by
\begin{equation}\label{rr:gridcount}
 \sum_{r=0}^{\min(2K,N)}\binom{N}{r}\bigl(2(2D^4+1)\bigr)^r
 \le \exp(CK\log D),
\end{equation}
since $N\le D$ and $D\ge2\mathrm e$. The factor $2$ inside the sum accounts for the allocation of each coordinate to either prefix or increment. The polynomial-size grid for $\lambda$ adds at most $C_h\log D$, since every coordinate standard deviation for an input in either unit ball is at most $\sqrt D$ under $M=1$, so $0\le\lambda\le C_hD^2$. The logarithm of the number of dyadic pairs $K,L\le D$ and choices of $Z$ is at most $C_h\log D$.

For a fixed grid triple the differences
\[
 V_i=\phi_\lambda(Z_i(w_1))-\phi_\lambda(Z_i(w_0))
\]
are independent across output coordinates, and $|V_i|\le |Z_i(e)|$. In the original orientation, let $V_i'$ be an independent copy of $V_i$. Jensen's inequality, symmetry of $V_i-V_i'$ and $\norm{V_i-V_i'}_{L_{2r}}\le 2\gamma_{2r}\sigma_i(e)$ give, for $u\in\R$,
\begin{align*}
 \int\exp\left(u\left(V_i-\int V_i\,d\prob\right)\right)d\prob
 &\le\int\exp(u(V_i-V_i'))\,d\prob\\
 &\le\sum_{r=0}^\infty\frac{(2u\sigma_i(e))^{2r}}{2^r r!}
 =\exp(2u^2\sigma_i(e)^2).
\end{align*}
Here $\gamma_{2r}^{2r}=(2r)!/(2^r r!)$ for $r\ge1$. By independence and $\sum_i\sigma_i(e)^2\le C_2^2\norm{e}_2^2$, the exponential Markov inequality, applied with both signs of $u$, yields
\begin{equation}\label{rr:originaltail}
 \prob\left(\left|\sum_i\left(V_i-\int V_i\,d\prob\right)\right|>C C_2\norm{e}_2\sqrt{z}\right)\le2e^{-z}\quad(z>0).
\end{equation}
If $C_2\norm{e}_2=0$, all the increments vanish.

In the transposed orientation put $b=C_2\norm{e}_\infty$. Every $Z_j(e)$ has standard deviation at most $b$. Restrict the union bound in this orientation to the deterministic triples for which
\[
 \sum_j\prob(E_j)^{1/2}\le4L,\qquad
 E_j=\{|Z_j(w_0)|>\lambda\}\cup\{|Z_j(w_1)|>\lambda\}.
\]
All triples arising from the threshold in \eqref{rr:activity} satisfy this condition. Here $V_j=0$ outside $E_j$. For each integer $m_0\ge2$, Cauchy--Schwarz, Gaussian moments and centering give
\begin{align*}
 &\sum_j\int_\Omega\left|V_j-\int_\Omega V_j\,d\prob\right|^{m_0}\,d\prob\\
 &\quad\le 2^{m_0}\sum_j\left(\int_\Omega|Z_j(e)|^{2m_0}\,d\prob\right)^{1/2}\prob(E_j)^{1/2}\\
 &\quad\le L\,m_0!\,(C b)^{m_0}.
\end{align*}
For the last inequality, one may use $\gamma_{2m_0}^{m_0}\le (2m_0)^{m_0/2}$ and $m_0!\ge(m_0/\mathrm e)^{m_0}$. Absorbing the numerical factor $4$ from the activity bound into $C^{m_0}$ leaves the prefactor $L$. In particular the centered independent summands satisfy the Bernstein moment condition for every integer $r\ge2$
\[
 \sum_j\int |V_j-\textstyle\int V_j\,d\prob|^r\,d\prob
 \le \tfrac12 r!\,\tau^2 K_0^{r-2},
 \qquad \tau^2=C_1 L b^2,\quad K_0=C_2' b,
\]
for sufficiently large absolute constants $C_1,C_2'>0$. Bernstein's inequality in its moment form \cite[Section 2.8]{BLM} gives
\begin{equation}\label{rr:bernstein}
 \prob\left(\left|\sum_j(V_j-\textstyle\int V_j\,d\prob)\right|>C b(\sqrt{Lu}+u)\right)\le2\exp(-u)\quad(u>0).
\end{equation} 

Take $u=C_h(z+1)K\log D$, with $z\ge1$. Including the threshold and level choices in \eqref{rr:gridcount} increases the constant in the exponential count. Choosing $C_h$ sufficiently large, the total exceptional probability is at most
\begin{equation}\label{rr:gridunion}
 C_h\sum_{K=2^j,\,j\ge0}\exp(-2(z+1)K\log D)
 \le C_h e^{-z},
\end{equation}
where the sum is bounded by the geometric series over all positive integers $K$. Also include the event $\max_{ij}|a_{ij}g_{ij}|\le C\sqrt{(z+1)\log D}$. With probability at least $1-C_he^{-z}$, the following estimates for the centered increments of $\Psi_\lambda$ hold simultaneously in absolute value, with respective upper bounds
\[
 \begin{gathered}
 C_h C_2\norm{e}_2\sqrt{(z+1)K\log D}\qquad\text{(original)},\\[2pt]
 C_h C_2\norm{e}_\infty\big(\sqrt{(z+1)LK\log D}+(z+1)K\log D\big)\\[-2pt]
 \text{(transposed)}.
 \end{gathered}
\]
In particular, these bounds apply when the input and its truncation level depend on the Gaussian matrix.

\smallskip
\noindent\emph{Summation over the dyadic shells.}
Order the absolute coordinates of $x\in B_p^n$ and $y\in B_v^m$ separately, breaking ties arbitrarily, and use dyadic rank shells $[K,2K)$. For a $y$ shell starting at $L$, collect all $x$ shells starting at $K\le L$. Its contribution in absolute value is at most $|y|_L^* T_L(G_Ax_{<2L})$, where the subscript retains the indicated input ranks (intersected with the available coordinates). For the complementary pairs use the transpose, an $x$ shell at $L$, and the $y$ shells at $K<L$. More explicitly, if $x^{[K]}$ and $y^{[L]}$ denote the two shell restrictions and $x_{<2L}=\sum_{K\le L}x^{[K]}$, $y_{<L}=\sum_{K<L}y^{[K]}$, then
\begin{equation}\label{rr:shellsplit}
 y^tG_Ax=\sum_L (y^{[L]})^tG_Ax_{<2L}
       +\sum_L (x^{[L]})^tG_A^ty_{<L}.
\end{equation}
All sums are over powers of two for which the corresponding restrictions are nonempty.

Use \eqref{rr:activity} with the full input $x$ or $y$, respectively. Let $w_0=0,w_1,\ldots,w_R$ be the successive rounded prefixes of the head in question, where $R$ is the number of its nonempty rank shells. With $\lambda$ chosen from the full input $w$, we have
\begin{align*}
 T_L(Z(w_R))&\le L\lambda+\psi_\lambda(w_R)\\
 &\quad+\sum_{a=1}^R\left[(\Psi_\lambda-\psi_\lambda)(w_a)
 -(\Psi_\lambda-\psi_\lambda)(w_{a-1})\right]\\
 &\le C_h H_L(\sigma(w))+LD^{-4}\\
 &\quad+\sum_{a=1}^R\left|(\Psi_\lambda-\psi_\lambda)(w_a)
 -(\Psi_\lambda-\psi_\lambda)(w_{a-1})\right|,
\end{align*}
by \eqref{rr:activity}. The difference between $T_L$ evaluated on the rounded and unrounded prefixes is at most $2L^2D^{-4}\max_{ij}|a_{ij}g_{ij}|$, since at most $2L$ input coordinates are rounded.

In the original orientation the rounded shell $e$ satisfies $\norm{e}_2\le K^{1/2-1/p}$. Multiplication by $|y|_L^*\le L^{-1/v}$ leaves a factor $K^{1-1/p}L^{-1/v}\le1$, since $K\le L$ and $1/p'\le1/v$. In the transposed orientation, $\norm{e}_\infty\le K^{-1/v}$ and $|x|_L^*\le L^{-1/p}$. The two factors left after multiplication are
\[
 K^{1/2-1/v}L^{1/2-1/p}\le1,\qquad
 K^{1-1/v} L^{-1/p}\le1 .
\]
These follow from $v>2$, $K\le L$ and $1/p+1/v\ge1$. The contribution of each pair of dyadic shells to the deviation is therefore at most $C_h(z+1)C_2\log D$. There are at most $C(\log D)^2$ such pairs. Writing $M_g=\max_{ij}|a_{ij}g_{ij}|$, the total rounding error is bounded by
\[
 C D^{-4}\sum_{L\le D,\,L\ {\rm dyadic}}(L+L^2 M_g)
 \le C(D^{-3}+D^{-2}M_g).
\]
On the event used in \eqref{rr:gridunion}, $M_g\le C\sqrt{(z+1)\log D}$ and $C_2\ge1$, so this error is bounded by $C(z+1)C_2\log D$. Estimate \eqref{rr:hardy} sums the original deterministic terms with $c=q$ and the transposed terms with $c=t$. Their total is bounded by $C_{p,q}(D_1(A)+D_2(A))$ by \eqref{d1} and \eqref{d2}. Consequently, for $z\ge1$,
\[
 \prob\bigl(\norm{G_A}_{p\to q}>C_{p,q}(D_1+D_2+(z+1)C_2\log^3 D)\bigr)
 \le C_{p,q}e^{-z}.
\]
Integration of this tail, using $\int_1^\infty e^{-z}\,dz<\infty$, gives the claimed estimate for $M=1$. Rescaling $A$ gives the general case.
\end{proof}

\subsection{Connected supports and the uniform weak variance}
Let $\Gamma(A)$ be the bipartite graph with row vertices $\{1,\ldots,m\}$, column vertices $\{1,\ldots,n\}$ (taken as a disjoint copy), and edges $(i,j)$ for which $a_{ij}\ne0$. Put
\begin{equation}\label{rr:degree}
 d=\max\left\{2,\max_i|\{j:a_{ij}\ne0\}|,\max_j|\{i:a_{ij}\ne0\}|\right\}.
\end{equation}
Write $L=\log(2\mathrm e mn)$, $M=\max_{i,j}|a_{ij}|$ and $C_2=\max_j\norm{(a_{ij})_i}_2$, where $\mathrm e=\exp(1)$.
\begin{lemma}\label{rr:connected}
For $1<p\le q<2$ and every $0<\varepsilon<1$,
\[
 \int_\Omega\norm{G_A}_{p\to q}\,d\prob\lesssim_{p,q,\varepsilon}
 D_1+D_2+b_*\sqrt L+
 d^{2\varepsilon}\sqrt{\log(2d)}\, C_2.
\]
\end{lemma}
\begin{proof}
The localization on connected supports in \cite[Section 3.3]{LS} is adapted here to the exponents $p$ and $v$. For $x\in B_p^n$, $y\in B_v^m$ and integers $h,e\ge0$, let
\[
 J_h=\{j:d^{-\varepsilon(h+1)}<|x_j|^p\le d^{-\varepsilon h}\},\qquad
 I_e=\{i:d^{-\varepsilon(e+1)}<|y_i|^v\le d^{-\varepsilon e}\}.
\]
If $j\in J_h$, $i\in I_e$ and $h-e\ge1+\lceil p/\varepsilon\rceil$, then
\[
 |x_j|^p\le d^{-\varepsilon h}\le d^{-p}d^{-\varepsilon(e+1)}\le d^{-p}|y_i|^v.
\]
As $1/p+1/v\ge1$ and $|y_i|\le1$, we have $|x_jy_i|\le d^{-1}|y_i|^v$. Let $E_1$ be the set of support edges satisfying this level inequality, and put $M_g=\max_{ij}|a_{ij}g_{ij}|$. By \eqref{rr:degree},
\[
 \sum_{(i,j)\in E_1}|a_{ij}g_{ij}x_jy_i|
 \le M_g d^{-1}\sum_i d|y_i|^v\le M_g.
\]
For the set $E_2$ of edges with $e-h\ge1+\lceil v/\varepsilon\rceil$, we have $|y_i|^v\le d^{-v}|x_j|^p$ and $|x_jy_i|\le d^{-1}|x_j|^p$. Hence
\[
 \sum_{(i,j)\in E_2}|a_{ij}g_{ij}x_jy_i|
 \le M_g d^{-1}\sum_j d|x_j|^p\le M_g.
\] For the remaining pairs of levels, the integer $h-e$ lies in an interval with at most $2+\lceil p/\varepsilon\rceil+\lceil v/\varepsilon\rceil$ integers.

Fix one such difference and take the connected components with at least one edge of the graphs on $J_h\sqcup I_e$ induced by the support. For a component with $k$ columns and $l$ rows, normalize the two restrictions in $\ell_p$ and $\ell_v$. Their coordinate moduli are at most $d^\varepsilon k^{-1/p}$ and $d^\varepsilon l^{-1/v}$, since the original coordinates lie in the respective levels. Its $k$ column vertices form a connected set in the graph joining columns at distance two in the original support. Its $l$ row vertices have the corresponding property. Let $\mathcal J_k$ and $\mathcal I_l$ be the families of connected vertex sets of cardinalities $k$ and $l$ in the column and row graphs just defined. Put
\begin{align*}
 \mathcal X_J&=\{x\in B_p^n:\supp x\subseteq J,\ \norm{x}_\infty\le d^\varepsilon k^{-1/p}\},\\
 \mathcal Y_I&=\{y\in B_v^m:\supp y\subseteq I,\ \norm{y}_\infty\le d^\varepsilon l^{-1/v}\},\\
 Y_{kl}&=\sup_{J\in\mathcal J_k,\ I\in\mathcal I_l}\ \sup_{x\in\mathcal X_J,\ y\in\mathcal Y_I}|y^tG_Ax|.
\end{align*}
The supremum is defined to be zero if either family is empty. Let $\mathcal C$ denote the family of components for the fixed difference, and write $J_C,I_C$ for the column and row sets of $C\in\mathcal C$. Put $\alpha_C=\norm{x_{J_C}}_p$ and $\beta_C=\norm{y_{I_C}}_v$. The sets are pairwise disjoint on each side, so
\begin{align*}
 \sum_{C\in\mathcal C}\alpha_C\beta_C
 &\le\left(\sum_C\alpha_C^p\right)^{1/p}
       \left(\sum_C\beta_C^{p'}\right)^{1/p'}\\
 &\le\left(\sum_C\alpha_C^p\right)^{1/p}
       \left(\sum_C\beta_C^v\right)^{1/v}\le1,
\end{align*}
since $p'\ge v$. Thus, pathwise,
\[
 \norm{G_A}_{p\to q}\le
 2\max_{ij}|a_{ij}g_{ij}|+C_{p,q,\varepsilon}\max_{k,l}Y_{kl}.
\]

The standard deviation of every scalar functional in $Y_{kl}$ is at most
\[
 \min(b_*,V_{kl}),\qquad
 V_{kl}=d^{2\varepsilon} k^{1/2-1/p}l^{-1/v}C_2.
\]
For the first bound use $\norm{x}_2\le1$ and $\norm{(a_{ij}y_i)_i}_2\le b_j\norm{y}_v$. For the second use both coordinate bounds and $\sum_{i\in I,j\in J}a_{ij}^2\le k C_2^2$.

If $k\le l$, use a half-net on each connected column support $J$ in the norm whose unit ball is
\[
 \{x:\supp x\subseteq J,\ \norm{x}_p\le1,\ \norm{x}_\infty\le d^\varepsilon k^{-1/p}\}.
\]
The norm on $\R^J$ is $\max\{\norm{x}_p,d^{-\varepsilon}k^{1/p}\norm{x}_\infty\}$. A half-net is a subset of the unit ball such that every point of the ball is within distance $1/2$ of a net point. The volume bound gives a half-net $\mathcal N_J$ with $|\mathcal N_J|\le5^k$. For a fixed $J$, put
\[
 M_J=\sup_{x\in\mathcal X_J,\ I\in\mathcal I_l,\ y\in\mathcal Y_I}|y^tG_Ax|.
\]
For $x\in\mathcal X_J$ choose $x_0\in\mathcal N_J$ with $2(x-x_0)\in\mathcal X_J$. Then
\[
 |y^tG_Ax|\le |y^tG_Ax_0|+\tfrac12M_J,
 \qquad M_J\le 2\max_{x_0\in\mathcal N_J}\sup_{I\in\mathcal I_l,\ y\in\mathcal Y_I}|y^tG_Ax_0|.
\] The number of supports is at most $nd^{4k}$, since a walk of $2(k-1)$ steps in the graph of degree at most $d^2$ traverses a spanning tree of each support. For fixed $x$ on the net, maximize over the entire admissible family of $y$. This is one Gaussian Lipschitz functional with mean at most $\gamma_qD_1$ by \eqref{rr:fixed} and Lipschitz constant at most $\min(b_*,V_{kl})$. When $l<k$, take the half-net on the row support and maximize over the columns. There are at most $m(5d^4)^l$ net points in this orientation and the mean bound is $\gamma_tD_2$. Applying \eqref{rr:conc} to these finite families and taking the union over $k,l$ gives
\begin{equation}\label{rr:connectedmax}
 \int_\Omega\max_{k,l}Y_{kl}\,d\prob\lesssim_{p,q}D_1+D_2+
 b_*\sqrt L+\sqrt{\log(2d)}\max_{k,l}\sqrt{\min(k,l)}V_{kl}.
\end{equation}
Indeed, the logarithm of the count for each pair is at most
\[
 CL+C\min(k,l)\log(2d).
\]
For $z\ge1$, take the concentration parameter to be $C(L+\min(k,l)\log(2d))+z$. The union bound over all size pairs and net points gives an exceptional probability at most $C_0e^{-z}$ after increasing the absolute constant $C$. Since the Lipschitz constant is at most $\min(b_*,V_{kl})$, the deviation outside this event is at most a constant times
\[
 b_*\sqrt{L+z}+V_{kl}\sqrt{\min(k,l)\log(2d)}.
\]
Integrating the tail gives \eqref{rr:connectedmax}.

For fixed $l$, the factor in $k$ satisfies
\[
 k^{1/2-1/p}\sqrt{\min(k,l)}
 =\begin{cases}k^{1-1/p},&k\le l,\\ l^{1/2}k^{1/2-1/p},&k\ge l.\end{cases}
\]
The first expression increases in $k$, whereas the second decreases, since $1<p<2$. Hence
\[
 \sqrt{\min(k,l)}V_{kl}
 \le d^{2\varepsilon} C_2\,l^{1/q-1/p}\le d^{2\varepsilon} C_2 .
\]
Since $|a_{ij}|\le b_*$, a union bound for the Gaussian tails gives
\[
 \prob\left(\max_{i,j}|a_{ij}g_{ij}|>u\right)\le 2mn\exp(-u^2/(2b_*^2))
\]
when $b_*>0$. Integrating, with a split at $b_*\sqrt{2\log(2\mathrm e mn)}$, bounds the expected largest entry by $C_0 b_*\sqrt L$. The case $b_*=0$ is immediate. The lemma follows.
\end{proof}

\begin{proposition}\label{rr:uniform}
For every $1<p\le q<2$ one has
\[
 \int_\Omega\norm{G_A}_{p\to q}\,d\prob\lesssim_{p,q}
 D_1(A)+D_2(A)+b_*\sqrt{\log(2\mathrm e mn)} .
\]
\end{proposition}
\begin{proof}
For a nonzero matrix normalize $D_1(A)+D_2(A)+b_*\sqrt L=1$ and put $\theta=1-q/2>0$. Then $M\le b_*\le L^{-1/2}$. Every restriction $A'$ obtained by erasing entries satisfies
\begin{equation}\label{rr:gain}
 C_2(A')\le D_1(A')^{q/2} M(A')^{\theta},
\end{equation}
because $\sum_i |a'_{ij}|^2\le M(A')^{2-q}\sum_i |a'_{ij}|^q$.

Set $K=4/\theta$. On the entries smaller than $L^{-K}$ in modulus, \eqref{rr:gain} gives $C_2\le L^{-K\theta}$. Lemma~\ref{rr:prelim} therefore bounds the expected norm of this part by a constant depending only on $p,q$.

On the complementary support, every column has $\ell_q$ norm at most $D_1\le1$ and every row has $\ell_t$ norm at most $D_2\le1$, the latter by testing $D_2$ on a coordinate vector. Its degree parameter satisfies $d\le 2L^{Kt}$, since $t>q$. By \eqref{rr:gain} one also has $C_2\le L^{-\theta/2}$. Choose $\varepsilon\in(0,1)$ so that $2\varepsilon Kt<\theta/4$. The degree term of Lemma~\ref{rr:connected} is bounded by
\[
 (2L^{Kt})^{2\varepsilon}\sqrt{\log(4L^{Kt})}\,L^{-\theta/2}
 \lesssim_{p,q}1 .
\]
All other terms decrease when entries are erased. Addition of the two estimates and homogeneity give the assertion.
\end{proof}

\subsection[The estimate with ordered column norms]{The estimate with ordered column norms}\label{sec:blocks}
We modify the construction in \cite[proof of Proposition 19]{LS} by including columns with the largest values of $b_j$ in the nested index sets.

Assume $A\ne0$ and normalize $S(A)=1$. Set $N_0=4$ and $N_k=N_{k-1}^2$ for $k\ge1$. Let $K_0$ be the least nonnegative integer with $N_{K_0}\ge\max(m,n)$. Append zero rows and columns to make $A$ square of order $N_{K_0}$. The quantities $D_1,D_2,B$ and the operator norm are unchanged. In this subsection $m=n=N_{K_0}$. Construct increasing row and column sets $I_k,J_k$ of sizes $N_k$, starting with arbitrary four-element sets $I_0,J_0$. Given $I_{k-1},J_{k-1}$, include in $I_k$ the rows of the $\lfloor N_{k-1}/4\rfloor$ largest entries in modulus of each column in $J_{k-1}$. Include in $J_k$ the columns of as many largest entries in each row of $I_{k-1}$, as well as the $N_{k-1}$ columns with largest $b_j$. The previous sets and all required additions use at most $2N_{k-1}+N_k/4\le N_k$ places on either side. Fill the unused places arbitrarily. For nonempty $I\subseteq\{1,\ldots,m\}$ and $J\subseteq\{1,\ldots,n\}$, write $G_{I,J}:\ell_p^J\to\ell_q^I$ for $(a_{ij}g_{ij})_{i\in I,j\in J}$. Put $I_{-1}=J_{-1}=\varnothing$, $P_k=I_k\setminus I_{k-1}$ and $Q_k=J_k\setminus J_{k-1}$. For $k\ge2$ one has
\begin{equation}\label{rr:thin}
 \begin{split}
 \max_{i\in P_k,\,j\in J_{k-2}} |a_{ij}|&\lesssim_q N_{k-2}^{-1/q},\qquad
 \max_{j\in Q_k,\,i\in I_{k-2}} |a_{ij}|\lesssim_p N_{k-2}^{-1/t},\\
 \max_{j\in Q_k} b_j&\le [\log(N_{k-2}+1)]^{-1/2}.
 \end{split}
\end{equation}
Indeed, at least $\lfloor N_{k-2}/4\rfloor$ entries in the indicated column or row have modulus no smaller than an entry outside the selected set. If $a$ is the modulus of such an entry, then $\lfloor N_{k-2}/4\rfloor a^q\le1$ in the column case and $\lfloor N_{k-2}/4\rfloor a^t\le1$ in the row case. Also, $J_{k-1}$ contains the $N_{k-2}$ columns with largest $b_j$; hence $b_j\le b_{N_{k-2}}^*\le[\log(N_{k-2}+1)]^{-1/2}$ for $j\in Q_k$.

Let $H$ be the least integer with $H\ge4$ and $2^{H-2}q/(2t)>1$. For a fixed offset $a=k-l$ with $|a|<H$, the rectangles $P_{l+a}\times Q_l$ have disjoint row and column sets. Put $M_l=\norm{G_{P_{l+a},Q_l}}_{p\to q}$ for $0\le l,l+a\le K_0$. For $x\in B_p^n$,
\[
 \left(\sum_l\norm{G_{P_{l+a},Q_l}x_{Q_l}}_q^q\right)^{1/q}
 \le\max_l M_l\left(\sum_l\norm{x_{Q_l}}_p^q\right)^{1/q}
 \le\max_l M_l.
\]
Testing on a vector supported in a single block gives the reverse inequality. Both $D_1$ and $D_2$ decrease upon restriction. Moreover $\log(2\mathrm e N_kN_l)\le C_H\log(N_{l-2}+1)$ for $l\ge2$. Proposition~\ref{rr:uniform} and \eqref{rr:thin} therefore bound the mean of every block norm. The first two column levels have bounded size in both dimensions for any such offset. The Gaussian Lipschitz constant of an operator norm in this range is at most its maximal column $\ell_s$ norm, by the first variance estimate in Lemma~\ref{rr:connected}. Consequently, for $l\ge2$ these constants are $O(2^{-l/2})$. Thus, if $C$ bounds the means of the block norms, \eqref{rr:conc} gives, for $z\ge1$,
\[
 \prob\bigl(\max_{l\ge2}M_l>C+z\bigr)
 \le \sum_{l\ge2}\exp(-c_H2^lz^2)\le C_H\exp(-c_Hz^2).
\]
Thus
\begin{align}
 \int_\Omega\max_l M_l\,d\prob
 &\le \sum_{l<2}\int_\Omega M_l\,d\prob
 +C+1+\int_1^\infty C_H\exp(-c_Hz^2)\,dz\\
 &\lesssim_{p,q,H}1.\label{rr:near}
\end{align}
Here the sum is over the admissible levels. An empty maximum is understood as zero.

Let $(X,\norm{\cdot})$ be a real normed space of dimension $d_0\ge1$, $F:\Omega\times X\to[0,\infty)$ a random seminorm and $f:X\to[0,\infty)$ a seminorm. Assume that $\gamma\ge0$, $\int_\Omega F(z)\,d\prob\le\gamma f(z)$ for $z\in X$, and $\sup_{\norm{z}\le1}f(z)\le1$. Assume also that for each $z\in B_X=\{z\in X:\norm{z}\le1\}$, $F(z)$ is a $\lambda$-Lipschitz function of a finite standard Gaussian vector. For $0<e\le1/2$ there is a random variable $E:\Omega\to[0,\infty)$ such that
\begin{equation}\label{rr:net}
 F(z)\le\gamma f(z)+E\norm{z}\quad\text{for every }z,\qquad
 \int_\Omega E\,d\prob\lesssim \lambda\sqrt{d_0\log(1+2/e)}+e\gamma .
\end{equation}
Indeed, let $\mathcal N$ be an $e$-net in the unit ball, with cardinality at most $(1+2/e)^{d_0}$, and set $Y=\max_{z_0\in\mathcal N}(F(z_0)-\gamma f(z_0))_+$. By \eqref{rr:conc},
\[
 \prob(Y>u)\le\min\{1,(1+2/e)^{d_0}\exp(-u^2/(2\lambda^2))\}.
\]
Integration gives $\int Y\,d\prob\le C\lambda\sqrt{d_0\log(1+2/e)}$ (the case $\lambda=0$ follows directly). Put $M_F=\sup_{\norm z\le1}F(z)$. For every $z$ in the unit ball, choose $z_0\in\mathcal N$ with $\norm{z-z_0}\le e$. The seminorm inequalities give
\begin{align*}
 F(z)&\le F(z_0)+eM_F\\
 &\le\gamma f(z)+Y+e(\gamma+M_F).
\end{align*}
It follows that $(1-e)M_F\le(1+e)\gamma+Y$, so $M_F\le3\gamma+2Y$. Substitution gives $F(z)\le\gamma f(z)+2Y+4e\gamma$ on the unit ball. Homogeneity proves \eqref{rr:net} with $E=2Y+4e\gamma$.

In the far upper part the disjoint output column sets $Q_k$ meet the input row prefixes $I_{k-H}$ for $k\ge H$. Define $F_k:\Omega\times\R^{I_{k-H}}\to[0,\infty)$ by $F_k(y)=\norm{G_{I_{k-H},Q_k}^{t}y}_t$ and $f_k:\R^{I_{k-H}}\to[0,\infty)$ by
\[
 f_k(y)=\left[\sum_{j\in Q_k}
 \left(\sum_{i\in I_{k-H}}a_{ij}^2y_i^2\right)^{t/2}\right]^{1/t}
 \qquad (y\in\R^{I_{k-H}}).
\]
By \eqref{rr:fixed}, $\int_\Omega F_k(y)\,d\prob\le\gamma_t f_k(y)$ and $\sup_{\norm{y}_v\le1}f_k(y)\le D_2(A)\le1$. For a rectangle with largest entry $M_0$, the Gaussian Lipschitz constant on the input unit ball is bounded by
\[
 \max_j\norm{(a_{ij})_i}_s\le
 M_0^{1-q/s}\max_j\norm{(a_{ij})_i}_q^{q/s}\le M_0^{q/2},
\]
with the columns restricted to that rectangle. Take $e=2^{-k-3}$. By \eqref{rr:thin} and \eqref{rr:net} the errors satisfy
\[
 \int_\Omega E_k\,d\prob\lesssim_{p,q}N_{k-2}^{-q/(2t)}\sqrt{N_{k-H}(k+1)}+2^{-k},
 \qquad \sum_{k\ge H}\int_\Omega E_k\,d\prob\lesssim_{p,q}1 ,
\]
because $N_{k-2}=N_{k-H}^{2^{H-2}}$ and $\delta=2^{H-2}q/(2t)-1/2>0$. In fact the first summand is $N_{k-H}^{-\delta}\sqrt{k+1}$, and $N_j=4^{2^j}$; the series converges independently of $K_0$. Let $U:\ell_p^n\to\ell_q^m$ be the random matrix obtained from $G_A$ by retaining only the rectangles $I_{k-H}\times Q_k$, $H\le k\le K_0$. For $y\in B_v^m$, \eqref{rr:net} and Minkowski's inequality give
\begin{align*}
 \norm{U^ty}_t
 &=\left(\sum_{k=H}^{K_0}F_k(y_{I_{k-H}})^t\right)^{1/t}\\
 &\le\gamma_t\left(\sum_{k=H}^{K_0}f_k(y_{I_{k-H}})^t\right)^{1/t}
       +\left(\sum_{k=H}^{K_0}E_k^t\norm{y_{I_{k-H}}}_v^t\right)^{1/t}\\
 &\le\gamma_t\norm{\widetilde u(y)}_t+\sum_{k=H}^{K_0}E_k.
\end{align*}
Thus
\begin{equation}\label{rr:farupper}
 \int_\Omega\norm{U}_{p\to q}\,d\prob
 =\int_\Omega\norm{U^t}_{v\to t}\,d\prob
 \le\gamma_tD_2(A)+\sum_{k=H}^{K_0}\int_\Omega E_k\,d\prob
 \lesssim_{p,q}1.
\end{equation}

Let $W:\ell_p^n\to\ell_q^m$ be the random matrix obtained from $G_A$ by retaining only the rectangles $P_k\times J_{k-H}$, $H\le k\le K_0$. Define $F_k^-:\Omega\times\R^{J_{k-H}}\to[0,\infty)$ and $f_k^-:\R^{J_{k-H}}\to[0,\infty)$ by
\begin{align*}
 F_k^-(x)&=\norm{G_{P_k,J_{k-H}}x}_q,\\
 f_k^-(x)&=\left[\sum_{i\in P_k}
       \left(\sum_{j\in J_{k-H}}a_{ij}^2x_j^2\right)^{q/2}\right]^{1/q}.
\end{align*}
By \eqref{rr:fixed},
\[
 \int_\Omega F_k^-(x)\,d\prob\le\gamma_q f_k^-(x),
 \qquad \sup_{x\in B_p^{J_{k-H}}}f_k^-(x)\le D_1(A)\le1.
\]
By \eqref{rr:thin} and the Lipschitz estimate used for $F_k$, the Lipschitz constant of $F_k^-(x)$, for $x\in B_p^{J_{k-H}}$, is at most
\[
 \left(\max_{i\in P_k,\,j\in J_{k-H}}|a_{ij}|\right)^{q/2}
 \lesssim_q N_{k-2}^{-1/2}.
\]
Applying \eqref{rr:net} with $e=2^{-k-3}$ gives a nonnegative random variable $E_k^-$ such that
\begin{align*}
 F_k^-(x)&\le\gamma_q f_k^-(x)+E_k^-\norm{x}_p,\\
 \int_\Omega E_k^-\,d\prob
 &\lesssim_{p,q}N_{k-2}^{-1/2}\sqrt{N_{k-H}(k+1)}+2^{-k}.
\end{align*}
Since $N_{k-2}=N_{k-H}^{2^{H-2}}$ and $H\ge4$,
\begin{align*}
 \sum_{k=H}^{K_0}\int_\Omega E_k^-\,d\prob
 &\lesssim_{p,q}\sum_{k\ge H}
 \left(N_{k-H}^{-(2^{H-2}-1)/2}\sqrt{k+1}+2^{-k}\right)\\
 &\lesssim_{p,q}1.
\end{align*}
For $x\in B_p^n$, the sets $P_k$ are disjoint, so
\begin{align*}
 \norm{Wx}_q
 &=\left(\sum_{k=H}^{K_0}F_k^-(x_{J_{k-H}})^q\right)^{1/q}\\
 &\le\gamma_q\left[\sum_{k=H}^{K_0}\sum_{i\in P_k}
 \left(\sum_{j\in J_{k-H}}a_{ij}^2x_j^2\right)^{q/2}\right]^{1/q}
 +\sum_{k=H}^{K_0}E_k^-\\
 &\le\gamma_q\norm{u(x)}_q+\sum_{k=H}^{K_0}E_k^-.
\end{align*}
Consequently,
\begin{equation}\label{rr:farlower}
 \int_\Omega\norm{W}_{p\to q}\,d\prob
 \le\gamma_qD_1(A)+\sum_{k=H}^{K_0}\int_\Omega E_k^-\,d\prob
 \lesssim_{p,q}1.
\end{equation}
For $|a|<H$, let $V_a:\ell_p^n\to\ell_q^m$ denote the random matrix obtained from $G_A$ by retaining the rectangles $P_{l+a}\times Q_l$ with $0\le l,l+a\le K_0$. The rectangles defining $U$, $W$ and the $V_a$ partition the matrix. By \eqref{rr:near}, \eqref{rr:farupper} and \eqref{rr:farlower},
\begin{align*}
 \int_\Omega\norm{G_A}_{p\to q}\,d\prob
 &\le\sum_{|a|<H}\int_\Omega\norm{V_a}_{p\to q}\,d\prob
   +\int_\Omega\norm{U}_{p\to q}\,d\prob
   +\int_\Omega\norm{W}_{p\to q}\,d\prob\\
 &\lesssim_{p,q}1.
\end{align*}
Removing the normalization $S(A)=1$ proves the upper bound in \eqref{remainingest} for $1<p\le q<2$.

If $2<p\le q<\infty$, then $1<q'\le p'<2$. By \eqref{duality} and the estimate just proved for $G_A^t$,
\begin{align*}
 \int_\Omega\norm{G_A}_{p\to q}\,d\prob
 &=\int_\Omega\norm{G_A^t}_{q'\to p'}\,d\prob\\
 &\lesssim_{p,q}D_1(A)+D_2(A)+L_{p,q}(A),
\end{align*}
since transposition interchanges $D_1$ and $D_2$ and the column exponent $2p'/(2-p')$ equals $2p/(p-2)$. The lower bound follows by the same duality.

\section{Proof of Theorem \ref{full} and Corollary \ref{infinite}}
\begin{proof}[Proof of Theorem \ref{full}]
For $q<p$ apply Theorem~\ref{main}; for $1<p\le q<2$ or $2<p\le q<\infty$ apply Theorem~\ref{mainremaining}. The range $1<p\le2\le q<\infty$ is covered by \cite[Theorem 2]{LS}; the remaining boundary cases follow from \cite[Propositions 1.8 and 1.10, Corollary 1.11 and (1.6)]{APSS}.
\end{proof}

\begin{proof}[Proof of Corollary \ref{infinite}]
By Theorem~\ref{full},
\[
 \sup_{N\ge1}\int_\Omega\norm{G_{A_N}}_{p\to q}\,d\prob
 \asymp_{p,q}\sup_{N\ge1}\bigl(D_1(A_N)+D_2(A_N)+L_{p,q}(A_N)\bigr).
\]
For each realization, the finite-section norms increase to the extended matrix norm of $G_A$. Indeed, if their supremum is $M<\infty$, then every truncated row has $\ell_{p'}$ norm at most $M$. Each row therefore acts by an absolutely convergent series on $\ell_p$, also when $p=\infty$. Passing to the limit on each finite set of output coordinates gives $\norm{G_A}_{p\to q}\le M$. The reverse bound follows by restriction. By taking countable dense subsets of each finite-dimensional pair of unit balls, the supremum of the finite-section norms is a countable supremum of centered Gaussian linear forms, including the negatives of these forms. The integrability of an almost surely finite Gaussian supremum (see \cite[proof of Proposition~1.2]{APSS}) shows that this supremum is integrable whenever it is finite almost surely. Monotone convergence now gives the equivalence with \eqref{infiniteest}.

Let $C=(c_{ij})_{i,j\ge1}$ be a nonnegative matrix and $C_N=(c_{ij})_{i,j\le N}$. Fix $a,b\in(0,\infty]$. Put $M=\sup_N\norm{C_N}_{a\to b}$. If $M<\infty$, then for every nonnegative $x\in\ell_a$ and $N\ge1$,
\[
 \norm{\left(\sum_{j=1}^N c_{ij}x_j\right)_{i=1}^N}_b
 \le M\norm{(x_j)_{j=1}^N}_a\le M\norm{x}_a.
\]
For $b<\infty$, monotone convergence gives
\[
 \sum_{i\ge1}\left(\sum_{j\ge1}c_{ij}x_j\right)^b
 =\lim_{N\to\infty}\sum_{i=1}^N\left(\sum_{j=1}^N c_{ij}x_j\right)^b
 \le M^b\norm{x}_a^b.
\]
For $b=\infty$,
\[
 \sup_{i\ge1}\sum_{j\ge1}c_{ij}x_j
 =\sup_{N\ge1}\max_{1\le i\le N}\sum_{j=1}^N c_{ij}x_j
 \le M\norm{x}_a.
\] In either case,
\[
 \norm{\left(\sum_{j\ge1}c_{ij}x_j\right)_{i\ge1}}_b
 \le M\norm{x}_a.
\]
Applying this to $|x|$ gives absolute convergence of the coordinate series for arbitrary $x\in\ell_a$ and $\norm{Cx}_b\le M\norm{x}_a$. Conversely, extending vectors by zero gives $\norm{C_N}_{a\to b}\le\norm{C}_{a\to b}$ whenever $C$ defines a bounded matrix operator. Hence, with infinite norms permitted,
\begin{equation}\label{inf:positive}
 \norm{C}_{a\to b}=\sup_{N\ge1}\norm{C_N}_{a\to b}.
\end{equation}
By \eqref{sq1}, \eqref{sq2} and \eqref{inf:positive},
\begin{align*}
 \sup_{N\ge1}D_1(A_N)^2&=\norm{A\circ A}_{p/2\to q/2},\\
 \sup_{N\ge1}D_2(A_N)^2&=\norm{(A\circ A)^t}_{q'/2\to p'/2}.
\end{align*}
For $p\le2\le q$, monotone convergence gives
\[
 \sup_{N\ge1}L_{p,q}(A_N)
 =\int_\Omega\sup_{i,j\ge1}|a_{ij}g_{ij}|\,d\prob.
\]

Suppose $p\le q<2$, and put $s=2q/(2-q)$ and
\[
 b_j^{(N)}=\left(\sum_{i=1}^N|a_{ij}|^s\right)^{1/s},
 \qquad b_j=\lim_{N\to\infty}b_j^{(N)}.
\]
For $r\ge1$ define
\[
 \eta_r=\sup_{\substack{J\subset\mathbb N\\ |J|=r}}\min_{j\in J}b_j.
\]
For $1\le r\le N$, let $(b^{(N)})_r^*$ be the $r$-th largest of $b_1^{(N)},\ldots,b_N^{(N)}$. Then
\[
 (b^{(N)})_r^*=\max_{\substack{J\subset\{1,\ldots,N\}\\|J|=r}}\min_{j\in J}b_j^{(N)}.
\]
Consequently,
\begin{align*}
 K:=\sup_{N\ge1}L_{p,q}(A_N)
 &=\sup_{r\ge1}\sqrt{\log(r+1)}\,
   \sup_{N\ge r}\max_{\substack{J\subset\{1,\ldots,N\}\\|J|=r}}
       \min_{j\in J}b_j^{(N)}\\
 &\le\sup_{r\ge1}\sqrt{\log(r+1)}\,\eta_r.
\end{align*}
Conversely, for any finite $J\subset\mathbb N$ with $|J|=r$ and $N\ge\max J$,
\[
 \sqrt{\log(r+1)}\min_{j\in J}b_j^{(N)}\le K.
\]
Letting $N\to\infty$ and then taking the supremum over $J$ and $r$ gives
\begin{equation}\label{inf:ordered}
 K=\sup_{r\ge1}\sqrt{\log(r+1)}\,\eta_r.
\end{equation}
If $K<\infty$, then $b_j\le K/\sqrt{\log 2}$ for every $j$. For each $\delta>0$, the existence of $r$ distinct indices with $b_j\ge\delta$ would imply $\delta\sqrt{\log(r+1)}\le K$ by \eqref{inf:ordered}. Thus
\[
 \big|\{j:b_j\ge\delta\}\big|\le \exp(K^2/\delta^2)-1<\infty.
\]
In particular, $b_j\to0$, and $\eta_r=b_r^*$ for every $r\ge1$. Conversely, if $b_j\to0$ and $\sup_r b_r^*\sqrt{\log(r+1)}<\infty$, then $K<\infty$ by \eqref{inf:ordered}.

If $2<p\le q$, apply \eqref{inf:ordered} to $A^t$ with exponents $(q',p')$: the corresponding column norms are
\[
 \left(\sum_{j\ge1}|a_{ij}|^{2p'/(2-p')}\right)^{(2-p')/(2p')}
 =\norm{(a_{ij})_{j\ge1}}_{2p/(p-2)}=d_i.
\]
At $p=\infty$, the exponent $2p/(p-2)$ is understood as $2$. For $q<p$ the term $L_{p,q}(A_N)$ vanishes for all $N$.
\end{proof}

\section*{Acknowledgments}
D. M. Pellegrino acknowledges support from CNPq through Grants No. 406457/2023-9 (CNPq/MCTI No. 10/2023); 403964/2024-5 (MCTI/CNPq No. 16/2024); 305807/2025-0. E. V. Teixeira gratefully acknowledges support from the Grayce B. Kerr Chair funds at Oklahoma State University. This research was conducted in part under the DARPA ExpMath project ``A Human-Centered Framework for AI-Mathematician Collaboration in Research-Level Mathematics'' (Agreement No. HR0011262E029), in which E. V. Teixeira serves as a co-principal investigator. He gratefully acknowledges DARPA's support and the collaborative environment fostered by the project. The views, opinions, and findings expressed here are those of the authors and should not be interpreted as representing the official views or policies of the Department of Defense or the U.S. Government.

\section*{Competing interests}
The authors declare that they have no competing interests.

\section*{Data availability}
No datasets were generated or analyzed in this theoretical study.

\section*{Declaration on the use of generative AI}
We began studying these conjectures in 2024 in connection with our research on the Kahane--Salem--Zygmund and Bennett inequalities. After exploring several approaches to the upper bounds, we used ChatGPT (OpenAI) to test estimates on examples and at endpoint exponents, and to develop parts of the proofs. ChatGPT contributed in particular to the use of a common truncation level for successive coordinate restrictions in the proof for $1<p\le q<2$. We also used ChatGPT for routine \LaTeX{} editing. We treated all AI-generated suggestions as provisional, checked the
arguments retained in the paper, and take full responsibility for the
final manuscript.

\end{document}